\documentclass[12pt]{article}

\usepackage{amssymb}
\usepackage{amsthm}
\usepackage{amsmath}
\usepackage{amsfonts}
\usepackage{hyperref}
\usepackage{color}
\usepackage[T1]{fontenc}
\usepackage[latin2]{inputenc}
\title{On some nonlinear extensions  of the Gagliardo-Nirenberg inequality with applications to nonlinear eigenvalue problems}
\author{ Agnieszka Ka\l{}amajska\thanks{The work is supported by
the  Polish
Ministry of Science grant no. N N201 397837 (years 2009-2012). }
 and Jan Peszek\thanks{The work is supported by EU FP6 Marie Curie RTN programme CODY}
 \\[2mm]
{\small \begin{tabular}{c}Faculty of Mathematics, Informatics, and Mechanics,
University of Warsaw, ul.
Banacha 2,\\ 02--097 Warszawa, Poland
\end{tabular}}
 \\[3mm]
\small email addresses: (A.K.) {\tt kalamajs@mimuw.edu.pl}, \small(J.P.) {\tt jp234556@students.mimuw.edu.pl}}

\date{}

\renewcommand{\it}{\sl}
\renewcommand{\em}{\sl}

plus 6pt minus 12pt
\newcommand{\barint}{
         \rule[.036in]{.12in}{.009in}\kern-.16in
          \displaystyle\int  }
\def\r{{\mathbb{R}}}
\def\n{{\mathbb{N}}}

\begin{document}

\newtheorem{theo}{\bf Theorem}[section]
\newtheorem{coro}{\bf Corollary}[section]
\newtheorem{lem}{\bf Lemma}[section]
\newtheorem{rem}{\bf Remark}[section]
\newtheorem{defi}{\bf Definition}[section]
\newtheorem{ex}{\bf Example}[section]
\newtheorem{fact}{\bf Fact}[section]
\newtheorem{prop}{\bf Proposition}[section]
\newtheorem{prob}{\bf Problem}[section]

\makeatletter \@addtoreset{equation}{section}
\renewcommand{\theequation}{\thesection.\arabic{equation}}
\makeatother

\newcommand{\ds}{\displaystyle}
\newcommand{\ts}{\textstyle}
\newcommand{\ol}{\overline}
\newcommand{\wt}{\widetilde}
\newcommand{\ck}{{\cal K}}
\newcommand{\ve}{\varepsilon}
\newcommand{\vp}{\varphi}
\newcommand{\pa}{\partial}
\newcommand{\rp}{\mathbb{R}_+}
\newcommand{\hh}{\tilde{h}}
\newcommand{\HH}{\tilde{H}}
\newcommand{\ct}{{\cal T}}

\maketitle

\smallskip
  {\small Key words and phrases: Gagliardo-Nirenberg  inequalities,
  interpolation inequalities, nonlinear eigenvalue problems\\

  }

  \begin{abstract}
  \noindent
  We derive inequality \[
  \int_{\r} |f^{'}(x)|^ph(f(x))dx \le \left(\sqrt{p-1}\right)^p \int_{\r}  \left( \sqrt{|f^{''}(x){\cal T}_h(f(x))| }\right)^ph(f(x))dx,\]
   where $f$ belongs locally to Sobolev space $W^{2,1}$ and $f^{'}$ has bounded support.
  Here $h(\cdot )$ is a given function and ${\cal T}_h(\cdot )$ is its given transform, it is independent of $p$.
   In case  when $h\equiv 1$ we retrieve the well known inequality:
  \(
  \int_{\r} |f^{'}(x)|^pdx \le \left(\sqrt{p-1}\right)^p \int_{\r} \left( \sqrt{|f^{''}(x)f(x)| }\right)^pdx.\)
  Our inequalities have  form similar to the classical second order Oppial inequalites. They also  extend certain class of inequalities due to Mazya,
  used to obtain  second order isoperimetric inequalities and capacitary estimates. We apply them to obtain new apriori estimates for nonlinear eigenvalue problems.
  %Orlicz extensions of the above inequality are also given.
\end{abstract}

\section{Introduction}
The purpose of this paper is two-fold. Our first goal is to derive certain generalization of the  classical interpolation inequality due to Gagliardo and Nirenberg  (\cite{g,n1}),  the second one is to present an application of the  inequality derived to  nonlinear elliptic eigenvalue problems in second order ODEs.

Let us start from explanation of our first approach.
It is well known that the following variant of the classical  Gagliardo-Nirenberg  inequality holds:
\begin{lem}\label{motyw}
For any $f\in C_0^\infty (\r)$, $p\ge 2$,    we have
\begin{eqnarray}
\int_\r |f^{'}|^pdx \le \left( \sqrt{p-1}\right)^p \int_\r {|ff^{''}|}^{\frac{p}{2}}dx.\label{ll2}
\end{eqnarray}
\end{lem}
\noindent
To get it, we just note that
\(
%\begin{eqnarray*}
I:=\int_\r \left( |f^{'}|^{p-2}f^{'}\right)\cdot f^{'}dx =- \int_\r \left( |f^{'}|^{p-2}f^{'}\right)^{'}\cdot fdx= -\int_\r g^{'}fdx,
%\end{eqnarray*}
\)
where $g^{'}=(p-1)|f^{'}|^{p-2}f^{''}$. Consequently
\(
I= -(p-1) \int_\r |f^{'}|^{p-2}(ff^{''})dx\leq
(p-1)\int_\r |f^{'}|^{p-2}|ff^{''}|dx.
\)
Now it suffices to apply H\"older's inequality and rearrange. We refer to \cite{kp1} for Orlicz generalizations of (\ref{ll2}).

\noindent
Having this simple computation in mind, we ask the question whether it is possible to obtain a more general inequality:
\begin{eqnarray}\label{nier120}
\int_{(a,b)}|f^{'}|^ph(f)dx\leq C_p\int_{(a,b)}\left|f^{''}{\cal T}_h(f)\right|^\frac{p}{2}h(f)dx,
\end{eqnarray}
where  $f\in W_{loc}^{2,1}((a,b))$ and  obeys some additional assumptions, $h$ is a given  function,
${\cal T}_h(\cdot )$ is a certain  transform of $h$  such that for $h\equiv 1$ we have $T_h(\lambda)=\lambda$, and constant $C_p$ does not depend on $f$.

\noindent
In this paper we give an affirmative answer to this question, constructing the appropriate operator ${\cal T}_h(\cdot)$ and proposing the  class of admissible functions $h$ and $f$.

\noindent
Our issue can be tracked back to Opial and  Mazy'a. Indeed, Opial obtained inequalities in the form:
\begin{equation}\label{opial11}
\int_a^b |yy^{'}|dx\le K \left( \int_a^b |y^{''}|^p dx\right)^{\frac{2}{p}},
\end{equation}
known as second order Opial inequalities, holding on compact interval $[a,b]$,
where
\[
y\in BC_0=\{ y\in W^{2,p}((a,b)): y(a)=y(b)=y^{'}(b)=0\},
\]
see e.g. see \cite{opial} and  e.g. \cite{bloom}, \cite{bbch}, \cite{fitz} for further issues and applications related to this inequalities.
 Inequality (\ref{opial11}) is similar to (\ref{nier120}), as the left hand side depends
on the functions  $y$ and $y^{'}$.

\noindent
Some other inequalities, having the form:
\begin{equation}\label{mazya}
\int_{{\rm supp} f^{'}} \left(\frac{ |f^{'}|}{f^{\frac{1}{2}}}\right)^pdx \le \left( \frac{p-1}{|1-\frac{1}{2} p|}\right)^{\frac{p}{2}}
 \int_\r |f^{''}|^{{p}}dx,
\end{equation}
holding for all smooth, nonnegative, compactly supported functions $f$,  can be found in Mazya's book \cite{ma}, Lemma 1, Section 8.2.1.
Those inequalities were the key arguments to obtain second order isoperimetric inequalities and capacitary estimates.
See also Proposition 1.9 in \cite{mazjakufner} for its weighted variant with $A_p$-Muckenhoupt weights and to Section 2 in the same article for its extensions
to  nonnegative functions of several variables.
 Let us mention that Mazya inequalities (\ref{mazya}) follow as a special case from our general result (\ref{nier120}), with the same constants.

\noindent
Our inequality (\ref{nier120}) implies inequality:
\[
\left( \int_\r |f^{'}|^{p}h(f) dx\right)^{\frac{2}{p}} \le (p-1) \left(  \int_\r |{\cal T}_h(f)|^qh(f)dx\right)^{\frac{1}{q}} \left( \int_\r |f^{''}|^rh(f)dx\right)^{\frac{1}{r}},
\]
where $q\ge \frac{p}{2}$, $\frac{2}{p}=\frac{1}{q}+\frac{1}{r}$. Substituting $h\equiv 1$ we obtain a generalization of the classical Gagliardo-Nirenberg inequality:
\[
\left( \int_\r |f^{'}|^{p} dx\right)^{\frac{2}{p}} \le C \left(  \int_\r |f|^qdx\right)^{\frac{1}{q}} \left( \int_\r |f^{''}|^rdx\right)^{\frac{1}{r}},
\] see Remark \ref{uwaga4} for details.

\noindent
Taking into account possible further applications we tried to obtain inequalities (\ref{nier120}) in the fullest possible generality. Sections \ref{three}, \ref{four} and \ref{five}
are devoted to the derivation of the inequalities under various constraints (Section \ref{three}: with regular weights, Section \ref{four}: with irregular weights
and nonnegative functions, Section \ref{five}: with irregular weights and  possibly sign changing functions).
Section \ref{six} closes  our analysis on inequalities by providing examples. To the best of our knowledge  inequalities (\ref{nier120}), except the just discussed cases, were not known earlier.

We are now to explain our second purpose, applications to ODEs.
Let us consider the following eigenvalue problem:
\begin{eqnarray}\label{eq120}
\left\{
\begin{array}{ccc}
f^{''}(x)&=&g(x)\tau(f(x))\ \hbox{\rm a.e. in}\ (a,b),\\
f\in {\cal R}
\end{array}
\right.
\end{eqnarray}
where $-\infty\le a<b\le +\infty$, $\tau: A\rightarrow \r$,  $A$ is some subinterval of $\r$ (finite or not), $g\in L^q(a,b)$ with $q\in [1,\infty ]$ and
 $f\in  W^{2,1}_{loc}((a,b))$, $f(x)\in A$ and
set ${\cal R}$ defines the boundary conditions.

\noindent
Such equations appear in many mathematical models, for example those describing problems arising
in the study of plasma physics, in determining the electrical potential in an isolated neutral atom, in the theory of vibrating string and of
 shallow membrane caps, in the theory of colloids, in the theory describing  flow and heat transfer over a stretching sheet,
or unsteady flow of a gas through a porous medium, in the models of catalytic theory, chemically reacting systems and adiabatic tubular reactors,
as well as  in the mathematical biology.
  We refer e.g. to books: \cite{ar,arw,bgm,oregan},   or an overview article \cite{hmv} for
 description of the models, existence theory and more references.
This type of equations are  helpful
 in functional analysis to describe best constants in Sobolev - Poincare - type inequalities, see e.g.
   \cite{bgm}, Chapter 5.

\noindent
Using inequalities (\ref{nier120}) we obtain new regularity results for solutions of (\ref{eq120}).
Roughly speaking, our approach  can be explained in the following way. If we succeed to find a function $h(\cdot)$ such that
\[
|g(x)|^q=\left|\frac{f^{''}(x)}{\tau(f(x))}\right|^q=|\ct_h(f(x))f^{''}(x)|^\frac{2q}{2}h(f(x)),
\]
and if we can use (\ref{nier120}), then we are able to deduce that $\int |f^{'}|^{2q}h(f)\leq C\|g\|_q^q$. Automatically it follows that $G(f)$ is $\lambda$-H\"older continuous, where $\lambda
=1-\frac{1}{2q}$ and  $G=G_\tau$ is such transform of $\tau$  that $|G^{'}|=|f^{'}|\cdot h(f)^{1/(2q)}$. Consequently this allows to deduce further regularity results, including the
analysis of an asymptotic behavior of solutions.
 Those results are presented in Section \ref{seven}.

\noindent
The discussion is provided by a detailed analysis within equations with homogeneous nonlinearity $\tau(\lambda)=\lambda^\alpha$ and positive functions. The model equation representing such approach would be
  Thomas and Fermi model found independently in 1927 to determine the electrical potential in an isolated neutral atom (see  \cite{ar}, page 121,
or \cite{arw}, page x,  and celebrated historical articles  \cite{fermi}, \cite{tomas}):
\begin{eqnarray*}
\left\{
\begin{array}{c}
y^{''}(t)= t^{\frac{1}{2}}y(t)^{\frac{3}{2}} , \ t\in (0,\infty),\\
y(0)=0, \lim_{t\to\infty} y(t)=0.
\end{array}
\right.
& {\rm or} &
\left\{
\begin{array}{c}
y^{''}(t)= t^{\frac{1}{2}}y(t)^{\frac{3}{2}} , \ t\in (0,1),\\
y(0)=0, y(1)=0.
\end{array}
\right.
\end{eqnarray*}
Another well known model is Emden-Fowler problem (see e.g. \cite{arw}, page x) which appears in various branches of fluid dynamics:
\begin{eqnarray*}
\left\{
\begin{array}{c}
y^{''}+\lambda q(x)y^{-\gamma}=0 , \ x\in (0,1),\gamma >0\\
y(0)=y(1)=0,
\end{array}
\right.
\end{eqnarray*}
Let us mention also the logistic equation present in  mathematical biology (\cite{afbr}, \cite{lixi}):
\[
u^{''} +af(x)u- b(x)u^p=0,\ x\in (0,\infty),\ p\ge 1,
\]
where we can put either: $a=0$ or $b\equiv 0$.

The presentation of   models related to homogeneous nonlinearity is provided in Section \ref{lustration}.

It is possible to generalize this approach in many directions: to deal with $n$-dimensional eigenvalue problems, also including systems,  analysis on Riemanian manifolds or Carnot-Carath\'{e}odory groups,  to obtain Orlicz variants of inequality (\ref{nier120}) and regularity results within Orlicz-Sobolev spaces, or to deal with weighted inequalities and its applications.
We hope that such generalizations will be an object of  future papers.

\section{Preliminaries and notation}
{\bf Notation.} Let $\Omega\subseteq \r^d$ be an open domain, $d\in \n$. We use  standard notation: $C_0^\infty (\Omega)$ to denote smooth compactly supported functions, $W^{m,p}(\Omega)$ and $W^{m,p}_{loc}(\Omega)$
to denote the  global and local Sobolev functions defined on $\Omega$, respectively.
%In our case will be mostly $\Omega =\r$.
If $A\subseteq \r$ and $f$ is defined on $A$,
by $f\chi_A$ we denote an extension of $f$ by zero outside set $A$.

\smallskip

\noindent
In the sequel we will use the following definition.

\begin{defi}\label{def1}
Let $-\infty\le A<B\le \infty$, $h:(A,B)\rightarrow [0,\infty)$ be a continuous function and let $H$ be a given locally absolutely continuous primitive of $h$. We define the transform of $h$:
\begin{eqnarray*}
{\cal T}_h(\lambda):=\left\{
\begin{array}{ccc}
\frac{H(\lambda)}{h(\lambda)} &  {\rm if}&  h(\lambda)\neq 0,\\
0 &{\rm if}& h(\lambda)=0,
\end{array}
\right.
\end{eqnarray*}
where $\lambda\in (A,B)$.
\end{defi}

\noindent
In the sequel we will use the following simple lemma.

\begin{lem}\label{lcomp}
If $f: [-R,R]\rightarrow \r$ is absolutely continuous with values in the interval $[\alpha,\beta]$ and $L: [\alpha,\beta] \rightarrow\r$ is Lipschitz, then
the function $(L\circ f)(x):=L(f(x))$ is absolutely continuous on $[-R,R]$.
\end{lem}

\section{Inequalities with regular weights $h(\cdot)$}\label{three}

\subsection{First approach. A model inequality}

Our first result reads as follows.

\begin{prop}\label{goal1}
Let $p\geq 2$, $h:\r\rightarrow[0,\infty)$ be a continuous function which is positive on $\r\setminus\{0\}$ and suppose that $H$ is an arbitrary chosen absolutely continuous primitive of $h$. Then for every $f\in {\cal R}$ where
\[
{\cal R}:= \{ f\in W_{loc}^{2,1}(\r) : \liminf_{r\to +\infty}|f^{'}(r)|^{p-1}|H(f(r))|= \liminf_{r\to -\infty}|f^{'}(r)|^{p-1}|H(f(r))|=0\},
\]
 we have
\begin{eqnarray}\label{nier1}
\int_\r|f^{'}(x)|^ph(f(x))dx\leq \left( \sqrt{p-1}\right)^p\int_\r\left|f^{''}(x){\cal T}_h(f(x))\right|^\frac{p}{2}h(f(x))dx,
\end{eqnarray}
where ${\cal T}_h$ is given by Definition \ref{def1}.
\end{prop}

\noindent
{\bf Proof of Proposition \ref{goal1}.}
Let us choose $r_n\to -\infty$ and $R_n\to +\infty$ such that $\lim_{n\to +\infty}|f^{'}(r_n)|^{p-1}|H(f(r_n))|= \lim_{n\to+\infty}|f^{'}(R_n)|^{p-1}|H(f(R_n))|=0$. Obviously we may assume that $I:=\int_\r|f^{'}|^ph(f)dx > 0$.  We have
\begin{eqnarray*}
I_n:=\int_{r_n}^{R_n}|f^{'}|^ph(f)\, dx
 = \int_{r_n}^{R_n}\left(|f^{'}|^{p-2}f^{'}\right)\cdot\left(h(f)f^{'}\right)\, dx.
\end{eqnarray*}
Function $f$ is  absolutely continuous on each interval $[r_n,R_n]$, thus there exist real $m=\inf_{x\in[r_n,R_n]}f(x)$ and $M=\sup_{x\in[r_n,R_n]}f(x)$. Furthermore, $H^{'}=h$ is  bounded on $[m,M]$ and thus the function $H$ is  Lipschitz on $[m,M]$.
 Lemma \ref{lcomp} implies that $H(f)$ is absolutely continuous on $[r_n,R_n]$. In particular its derivative can be computed almost everywhere and equals $h(f(x))f^{'}(x)$ a.e..
By almost the same argument, when $p\ge 2$ the function  $\Phi_p(\lambda)=|\lambda|^{p-2}\lambda$ is   Lipschitz on every interval $[-L,L]$ where $L>0$ and $f^{'}$ is locally absolutely continuous. Therefore  $\Phi_p(f^{'})=|f^{'}|^{p-2}f^{'}$ is absolutely continuous on $[r_n,R_n]$. We integrate by parts to obtain
\begin{eqnarray*}
I_n=-(p-1)\int_{r_n}^{R_n}\left(|f^{'}|^{p-2}f^{''}\right)H(f)\, dx+\theta_n \leq (p-1)\int_{r_n}^{R_n}\left(|f^{'}|^{p-2}|f^{''}|\right)|H(f)|\, dx +\theta_n,
\end{eqnarray*}
where
\[
 \theta_n = |f^{'}(R_n)|^{p-2}f^{'}(R_n) H(f(R_n)) - |f^{'}(r_n)|^{p-2}f^{'}(r_n) H(f(r_n)).
\]
Function $h(f(x))$ is positive a.e. on $\{x:f^{'}(x)\neq 0\}$. Indeed, points $x$ where $h(f(x))=0$  are contained in the set
$\{ x: f(x)=0\}$. On that set $f^{'}=0$ a.e..
 Thus $\frac{h(f(x))}{h(f(x))}=1$ almost everywhere on $\{x:f^{'}(x)\neq 0\}$ and
\begin{eqnarray*}
I_n\leq(p-1)\int_{(r_n,R_n)\cap \{ x: f^{'}(x)\neq 0\}}\left(|f^{'}|^{p-2}h(f)^\frac{p-2}{p}\right) \left(\frac{|f^{''}H(f)|}{h(f)^\frac{p-2}{p}}\right)dx +\theta_n.
\end{eqnarray*}
We apply H\"older inequality with $q=\frac{p}{p-2}$, $q^{'}=\frac{p}{2}$ to get
\begin{eqnarray}\label{potem}
I_n\leq (p-1)\left(\int_{r_n}^{R_n}|f^{'}|^p h(f) dx\right)^\frac{p-2}{p}\left(\int_{ (r_n,R_n)\cap \{ x: f^{'}(x)\neq 0\}}\frac{|f^{''}H(f)|^\frac{p}{2}}{h(f)^\frac{p-2}{2}} dx\right)^\frac{2}{p} +\theta_n=\nonumber\\
=(p-1)I_n^\frac{p-2}{p}\left(\int_{r_n}^{R_n}|f^{''}{\cal T}_h(f)|^\frac{p}{2}h(f)dx\right)^\frac{2}{p}+\theta_n.
\end{eqnarray}
Since we assumed that $I>0$ there exists some $n_0$ such that for $n>n_0$ we have $\infty>I_n\geq I_{n_0}>0$ and thus the sequence $\frac{\theta_n}{I_n^{(p-2)/p}}$ converges to $q\in [-\infty, 0]$. Therefore now it suffices to rearrange inequality (\ref{potem}) and converge with $n$ to infinity.\hfill
$\Box$

\begin{rem}\rm
Function $f\in W^{2,1}_{loc}(\r)$ belongs to ${\cal R}$ in the following situations:
 \begin{enumerate}
 \item when  $f^{'}$ has bounded support;
 \item
when
$
\int_\r |f^{'}(x)|^{p-1}|H(f(x))|\, dx <\infty$;
\item when  integrals:  $\int_\r |f^{'}|^ph(f)dx$ and $\int_\r |{\cal T}_h(f)|^ph(f)dx$ are finite.
\end{enumerate}
Indeed, first observation is obvious. Second one follows from the fact that\\
$\liminf_{r\to-\infty} |f^{'}(r)|^{p-1}|H(f(r))|= \liminf_{r\to-\infty} |f^{'}(r)|^{p-1}|H(f(r))|=0$, as otherwise integral $
\int_\r |f^{'}(x)|^{p-1}|H(f(x))|\, dx $ could not be finite.\\
Let us explain the last observation. As
\(
|f^{'}|^{p-1}|H(f)|=|f^{'}|^{p-1}|{\cal T}_h(f)|h(f),
\)
we get from H\"older inequality
\[
\int_\r |f^{'}|^{p-1}|H(f)| dx \le \left(\int_\r |f^{'}|^{p}h(f)dx   \right)^{1-\frac{1}{p}} \left( \int_\r|{\cal T}_h(f)|^{p}h(f)dx   \right)^{\frac{1}{p}}<\infty
\]
and we  apply second observation.
\end{rem}

\begin{rem}\label{uwaga4}\rm
Suppose that $f$ satisfies the assumptions of Proposition \ref{goal1}. Applying H\"older inequality with $\theta_1=(2q)/p> 1$ and $\theta_2=\theta_1/(\theta_1-1)$ to the right hand side of  (\ref{nier1}) and treating $h(f)dx$ as a measure,  we obtain inequality:
\[
\left( \int_\r |f^{'}|^{p}h(f) dx \right)^{\frac{2}{p}} \le (p-1)\left(  \int_\r |{\cal T}_h(f)|^qh(f)dx\right)^{\frac{1}{q}} \left( \int_\r |f^{''}|^rh(f)dx\right)^{\frac{1}{r}},
\]
where $r$ is such that $\frac{2}{p}=\frac{1}{q}+\frac{1}{r}$.
This is the nonlinear variant of classical Gagliardo-Nirenberg inequality:
\[
\left( \int_\r |f^{'}|^{p} dx\right)^{\frac{2}{p}} \le C\left(  \int_\r |f|^qdx\right)^{\frac{1}{q}} \left( \int_\r |f^{''}|^rdx\right)^{\frac{1}{r}}
\]
as ${\cal T}_{h\equiv 1}(f)= f$.
\end{rem}

\subsection{The generalization}

We are now to present the variant of Proposition \ref{goal1} which holds under minimal assumptions, provided that the function $h(\cdot)$ is continuous.
We assume that $A\le f\le B$ (possibly $A=-\infty, B=+\infty$) and we relax the assumption: $h(\lambda)>0$ for $\lambda\neq 0$.
Also, our function $f$ can be defined on an arbitrary interval $(a,b)$ and the admitted class of functions even in case  $(a,b)=\r$ can be wider.
Our result reads as follows.

\begin{prop}\label{goal2}
Let $p\geq 2$, $-\infty\le A<B\le+\infty$, $h:[A,B]\rightarrow [0,\infty)$ be such a continuous function that set $h^{-1}(0)$ is at most countable.
Suppose that $H$ is an arbitrary chosen locally absolutely continuous primitive of $h$, and let  ${\cal T}_h$ be given by Definition \ref{def1}. Further, let $-\infty\le a<b\le +\infty$ be the given numbers.
 Then we have
\begin{description}
\item[(i)] For every $f\in W^{2,1}_{loc}(a,b)$ such that $A\le f\le B$ and every subinterval $[r,R]\subseteq (a,b)$
\begin{equation}\label{nier50}
\int_r^R|f^{'}(x)|^ph(f(x))dx\leq \left( \sqrt{p-1}\right)^p\int_r^R\left|f^{''}(x){\cal T}_h(f(x))\right|^\frac{p}{2}h(f(x))dx + \theta(r,R),
\end{equation}
where $\theta(r,R)=  |f^{'}(R)|^{p-2}f^{'}(R) H(f(R)) - |f^{'}(r)|^{p-2}f^{'}(r) H(f(r))$;
\item[(ii)]
For every $f\in \widetilde{{\cal R}}_{[A,B]}(a,b)$
\begin{eqnarray}\label{nier12}
\int_a^b|f^{'}(x)|^ph(f(x))dx\leq \left( \sqrt{p-1}\right)^p\int_a^b\left|f^{''}(x){\cal T}_h(f(x))\right|^\frac{p}{2}h(f(x))dx,
\end{eqnarray}
where
{\small
\begin{eqnarray*}
\widetilde{{\cal R}}_{[A,B]}(a,b)&:=&
\{ f\in W_{loc}^{2,1}(a,b): A\le f\le B \ {\rm and}\\
&& \left( \liminf_{R\nearrow b}|f^{'}(R)|^{p-2}f^{'}(R)H(f(R))- \limsup_{r\searrow a}|f^{'}(r)|^{p-2}f^{'}(r)H(f(r))\right)\le 0\}.
\end{eqnarray*}
}

\end{description}
\end{prop}

\noindent
{\bf Proof.} The proof follows by an easy modification of the proof of Proposition \ref{goal1}. As set $ h^{-1}(0)$ is at most countable, set
${\cal B}=\{ x: h(f(x))=0\}$ is at most countable number of level sets of $f$. On each level set of $f$ we have $f^{'}=0$ almost everywhere.
Therefore $f^{'}=0$ almost everywhere on set ${\cal B}$.  Following the same lines as in the proof of Proposition \ref{goal1} we get inequality (\ref{nier50}),
 holding on every subinterval $[r,R]\subseteq (a,b)$.
To get (\ref{nier12}) suffices to chose the right sequence $r_n\searrow a$, $R_n\nearrow b$ such that $\lim_{n\to\infty}\theta (r_n,R_n)\le 0$, apply  inequality
 (\ref{nier50}) with $r=r_n,R=R_n$ and let $n$ converge to infinity.
\hfill$\Box$

\section{Inequalities dealing with less regular weights $h(\cdot)$ and nonnegative functions}\label{four}

Our next goal is to admit the more general class of functions $h$. We will deal with functions $h$ which may not be defined at zero and
the admitted class of functions consists of nonnegative functions only.

\noindent
We start with the following lemma which is a simple consequence of Proposition \ref{goal2}.

\begin{lem}\label{pierwszylemat}
 Let $p\geq 2$, $h: (0,\infty)\rightarrow (0,\infty)$ be a continuous function, $H$ and ${\cal T}_h$ be as in Definition \ref{def1},
 $-\infty \le a<b\le+\infty$, $\eta >0$.
 Then for every $f\in W^{2,1}_{loc}((a,b))$ such that $f\ge \eta$  and for every $r,R$ such that $a<r<R<b$, we have
\begin{eqnarray}\label{nier1111}
\int_r^R|f^{'}(x)|^ph(f(x))dx\leq \left( \sqrt{p-1}\right)^p\int_r^R\left|f^{''}(x){\cal T}_h(f(x))\right|^\frac{p}{2}h(f(x))dx +\theta(r,R),
\end{eqnarray}
where  \(
\theta (r,R) := |f^{'}(R)|^{p-2}f^{'}(R) H(f(R)) - |f^{'}(r)|^{p-2}f^{'}(r) H(f(r)).\)
\end{lem}

As direct consequence of the above Lemma one obtains the following result.

\begin{prop}\label{winosekdodatk}
Let $p\geq 2$, $h: (0,\infty)\rightarrow (0,\infty)$ be a continuous function, $H,{\cal T}_h: (0,\infty)\rightarrow\r$ be as in Definition \ref{def1},
 $-\infty \le a<b\le+\infty$.
 Then for every $f\in W^{2,1}_{loc}((a,b))$ such that $f>0$ on $(a,b)$, we have
\[
\int_a^b|f^{'}(x)|^ph(f(x))dx\leq \left( \sqrt{p-1}\right)^p\int_a^b\left|f^{''}(x){\cal T}_h(f(x))\right|^\frac{p}{2}h(f(x))dx +L_{-}(a,b)(f),
\]
where  \[
L_{-}(a,b)(f) := \left( \liminf_{R\nearrow b}|f^{'}(R)|^{p-2}f^{'}(R)H(f(R))- \limsup_{r\searrow a}|f^{'}(r)|^{p-2}f^{'}(r)H(f(r))\right). \]
In particular if $f\in \tilde{\cal R}_{>0}(a,b)$ where
\begin{eqnarray*}
 \tilde{\cal R}_{> 0}(a,b)&:=& \{ f\in W_{loc}^{2,1}((a,b)), f>  0\ {\rm on}\ (a,b) :\\
 &&
  \left( \liminf_{R\nearrow b}|f^{'}(R)|^{p-2}f^{'}(R){H}(f(R))
 -  \limsup_{r\searrow a}|f^{'}(r)|^{p-2}f^{'}(r){H}(f(r))\right) \le 0
\},
\end{eqnarray*}
 we have
\[
\int_a^b|f^{'}(x)|^ph(f(x))dx\leq \left( \sqrt{p-1}\right)^p\int_a^b\left|f^{''}(x){\cal T}_h(f(x))\right|^\frac{p}{2}h(f(x))dx.
\]
\end{prop}

\noindent
Some other additional assumptions allow to relax an assumption $f>0$ to the weaker one that $f$ is nonnegative almost everywhere. First result in this direction reads as follows.

\begin{prop}\label{goracs}
Let $p\geq 2$, $-\infty \le a<b\le+\infty$, $h:(0,\infty)\rightarrow(0,\infty)$ be a continuous function, $H,{\cal T}_h:(0,\infty)\rightarrow\r$ be as in Definition \ref{def1}.
Assume further that
 function $H$ extends to the continuous function $\tilde{H}: [0,\infty)\rightarrow\r$.
Then for every function $f\in \tilde{\cal R}^0_{\ge 0}(a,b)$ where
{\small
\begin{eqnarray*}
\tilde{\cal R}^0_{\ge 0}(a,b)&:=& \{ f\in W_{loc}^{2,1}((a,b)), f\ge  0\  a.e. \ {\rm on}\ (a,b)  :
%\label{er1}
\\ && \liminf_{R\nearrow b}|f^{'}(R)|^{p-2}f^{'}(R)\tilde{H}(f(R))\le 0
  \ {\rm and}\  \limsup_{r\searrow a}|f^{'}(r)|^{p-2}f^{'}(r)\tilde{H}(f(r)) \ge 0
\},
\end{eqnarray*}
}
we have
\begin{equation}\label{trzydwa}
\int_{(a,b)\cap \{ x:f(x)>0\}  }|f^{'}(x)|^ph(f(x))dx\leq \left( \sqrt{p-1}\right)^p\int_{(a,b)\cap \{ x:f(x)>0\}  }\left|f^{''}(x){\cal T}_h(f(x))\right|^\frac{p}{2}h(f(x))dx.
\end{equation}
\end{prop}

\begin{rem}\rm Let $h$ be as in Proposition \ref{goracs}. As $H$ is locally absolutely continuous on $(0,\infty)$, we have for every $y>0$: $H(x)-H(y)=\int_y^x h(\tau)d\tau$. After letting
$y\to 0$  we get that $H(x)=\int_0^xh(\tau)d\tau +\tilde{H}(0)$. In particular $h$ is integrable in  every neighborhood of $0$.
The converse is also true: if $h$ is integrable in some neighborhood of zero, then $H$ extends to absolutely continuous
function defined on $[0,\infty)$.
\end{rem}

\noindent
{\bf Proof of Proposition \ref{goracs}.}
 We can assume that $f$ has some zeroes on $(a,b)$, as otherwise the result follows directly from Proposition \ref{winosekdodatk}. Furthermore, since the integrands in (\ref{trzydwa}) are zero a.e. on the set $\{x\in (a,b): f^{'}(x)=0\}$,  we may consider only  integrals with respect to
 measure $\mu$ which is
  Lebesgue measure restricted to the set $\{x\in (a,b): f^{'}(x)\neq 0\}$.
By our assumptions set ${\cal S}:= \{ x\in (a,b): f(x)>0\}$ is of full measure $\mu$  and is a sum of disjoint open intervals $I$ which have one of the form:\\
(a) $I=(a,R)$ where $a<R<b$ and $f(R)=0$,\\
(b) $I=(r,R)$ where $a<r<R<b$ and $f(r)=f(R)=0$,\\
(c) $I=(r,b)$ where $a<r<b$ and $f(r)=0$.\\

\noindent
Proposition \ref{winosekdodatk} and the fact that the nonnegative function $f$ can have only double zeroes,   implies that
if $I$ is of the form (a), we have
{\small
\begin{equation}\label{sobota3}
\int_I|f^{'}(x)|^ph(f(x))dx\leq \left( \sqrt{p-1}\right)^p\int_I\left|f^{''}(x){\cal T}_h(f(x))\right|^\frac{p}{2}h(f(x))dx
- \limsup_{r\searrow a}|f^{'}(r)|^{p-2}f^{'}(r)H(f(r)).  \end{equation}
}
If $I$ is of the form (b), we have
\begin{equation}\label{sobota1}
\int_I|f^{'}(x)|^ph(f(x))dx\leq \left( \sqrt{p-1}\right)^p\int_I\left|f^{''}(x){\cal T}_h(f(x))\right|^\frac{p}{2}h(f(x))dx,
 \end{equation}
while if $I$ is of the form (c), we have
{\small
\begin{equation}\label{sobota2}
\int_I|f^{'}(x)|^ph(f(x))dx\leq \left( \sqrt{p-1}\right)^p\int_I\left|f^{''}(x){\cal T}_h(f(x))\right|^\frac{p}{2}h(f(x))dx +
 \liminf_{R\nearrow b}|f^{'}(R)|^{p-2}f^{'}(R)H(f(R)).\end{equation}
 }
Moreover, we have one of three possible situations:\\
(A) $a$ and $b$ are accumulation points of zeroes of $f$;\\
(B) exactly one of points $a$ or $b$ is an accumulation point of zeroes of $f$;\\
(C) $f$ is strictly positive in a neighborhood of $a$ and $b$. \\
In case (A) set ${\cal S}$ is a sum of intervals of the form (b), therefore we can sum up all inequalities which have the form (\ref{sobota1}) to obtain
{\small
\[
\int_{ (a,b)\cap \{x:f(x)>0\} }|f^{'}(x)|^ph(f(x))dx\leq \left( \sqrt{p-1}\right)^p\int_{ (a,b)\cap \{x: f(x)>0\} }\left|f^{''}(x){\cal T}_h(f(x))\right|^\frac{p}{2}h(f(x))dx.
 \]
 }
 In case (B), say when $a$ is an accumulation point of zeroes of $f$, we get from (\ref{sobota1}) and (\ref{sobota2}):
 {\small
 \begin{eqnarray*}
 &~&\int_{ (a,b)\cap \{x: f(x)>0\} }|f^{'}(x)|^ph(f(x))dx\leq \\ &~&\left( \sqrt{p-1}\right)^p\int_{ (a,b)\cap \{x: f(x)>0\} }\left|f^{''}(x){\cal T}_h(f(x))\right|^\frac{p}{2}h(f(x))dx +
\liminf_{R\nearrow b}|f^{'}(R)|^{p-2}f^{'}(R)\tilde{H}(f(R)).
 \end{eqnarray*}
 }
 In second case in (B) when $b$ is an accumulation point of zeroes of $f$, we get from (\ref{sobota1}) and and (\ref{sobota3}):
 {\small
\begin{eqnarray*}
 &~&
 \int_{ (a,b)\cap \{x:f(x)>0\} }|f^{'}(x)|^ph(f(x))dx\leq \\&~&\left( \sqrt{p-1}\right)^p\int_{ (a,b)\cap \{x:f(x)>0\} }\left|f^{''}(x){\cal T}_h(f(x))\right|^\frac{p}{2}h(f(x))dx -
 \limsup_{r\searrow a}|f^{'}(r)|^{p-2}f^{'}(r)\tilde{H}(f(r)).
 \end{eqnarray*}
 }
 Therefore in case (B) assertion follows.
 The case (C) follows by easy modifications of our previous arguments. \hfill$\Box$

Our next proposition applies to the case when $h$ is not necessarily integrable in a neighborhood of $0$.
It requires some additional assumptions on function $h(\cdot )$.

\begin{prop}\label{gorac}
Let $p\geq 2$, $-\infty \le a<b\le+\infty$, $h:(0,\infty)\rightarrow(0,\infty)$ be a continuous function, $H,{\cal T}_h:(0,\infty)\rightarrow\r$ be as in Definition \ref{def1}, $\tilde{H}(x)$ be an extension (not necessarily continuous) of $H$ up to $0$, $\tilde{H}(0)=0$.
Assume further that  both functions $h$ and $|{\cal T}_h|^\frac{p}{2}h$ are
 either bounded or nonincreasing in some neighborhood of $0$.
Then for every nonnegative function $f\in \tilde{\cal R}_{\ge 0}(a,b)$ where
{\small
\begin{eqnarray*}
\tilde{\cal R}_{\ge 0}(a,b)&:=& \{ f\in W_{loc}^{2,1}((a,b)), f\ge  0\  a.e. : \\
&&\left( \liminf_{R\nearrow b}|f^{'}(R)|^{p-2}f^{'}(R)\tilde{H}(f(R))- \limsup_{r\searrow a}|f^{'}(r)|^{p-2}f^{'}(r)\tilde{H}(f(r))\right)\le 0
\},
\end{eqnarray*}
}
we have
\[
\int_{(a,b)\cap \{x:f(x)>0\} }|f^{'}(x)|^ph(f(x))dx\leq \left( \sqrt{p-1}\right)^p\int_{(a,b)\cap \{ x:f(x)>0\} }\left|f^{''}(x){\cal T}_h(f(x))\right|^\frac{p}{2}h(f(x))dx.
\]
\end{prop}

\begin{rem}\rm Obviously, set
\begin{eqnarray*}
{\cal R}_{\ge 0}^0(a,b)&=& \{ f\in W_{loc}^{2,1}((a,b)), f\ge  0\ a.e. :
\liminf_{R\nearrow b}|f^{'}(R)|^{p-2}f^{'}(R)\tilde{H}(f(R))\le 0\ {\rm and}\\ && \limsup_{r\searrow a}|f^{'}(r)|^{p-2}f^{'}(r)\tilde{H}(f(r))\ge 0\}\end{eqnarray*}
is contained in $\tilde{\cal R}_{\ge 0}(a,b)$.
\end{rem}

\noindent
{\bf Proof of Proposition \ref{gorac}.}\\
We consider two cases:\\
{\sc Case A:} Both functions $h$ and $|{\cal T}_h|^\frac{p}{2}h$ are nonincreasing in some neighborhood of zero;\\
{\sc Case B:} One of functions $h$ or $|{\cal T}_h|^\frac{p}{2}h$ is bounded in some neighborhood of zero.\\

\smallskip
\noindent
{\sc Proof of Case A.}
We find some $\epsilon >0$ such that
 $h$ and $|{\cal T}_h|^\frac{p}{2}h$ are nonincreasing on $(0,\epsilon]$, denote $f_\eta(x):=f(x)+\eta$ where $0<\eta<\frac{\epsilon}{2}$, and
 chose $a<r<R<b$. Then we use Lemma \ref{pierwszylemat} and note that $f_\eta^{'}=f^{'}$,  to get
 \begin{eqnarray*}
  \int_r^R|f^{'}(x)|^ph(f_{\eta}(x))dx&\leq& \left( \sqrt{p-1}\right)^p\int_r^R\left|f^{''}(x){\cal T}_h(f_{\eta}(x))\right|^\frac{p}{2}h(f_{\eta}(x))dx +\theta_{\eta}(r,R),\\
  &=:& B_\eta (r,R) + \theta_\eta (r,R),\ {\rm where}\  \\
 \theta_\eta (r,R) &:=& |f^{'}(R)|^{p-2}f^{'}(R) H(f_{\eta}(R)) - |f^{'}(r)|^{p-2}f^{'}(r) H(f_{\eta}(r)).
 \end{eqnarray*}
 We will let $\eta$ converge to zero.
  We have
\begin{eqnarray}\label{tez1}
A_{\eta}(r,R):=\int_{\{x\in (r,R): f^{'}(x)\neq 0\}}|f^{'}|^ph(f_\eta)dx=\int_{E_\epsilon}|f^{'}|^ph(f_\eta)dx+ \int_{F_\epsilon}|f^{'}|^ph(f_\eta)dx,
\end{eqnarray}
where $E_\epsilon=\{x\in (r,R):0<f(x)<\frac{\epsilon}{2}\}$ and $F_\epsilon =\{x\in (r,R): f(x)\ge \frac{\epsilon}{2}\}$
(note that when $f(x)=0$ we necessarily have $f^{'}(x)=0$ as nonnegative functions can have double zeroes only).
 The former integrand increases as $\eta\rightarrow 0$. Furthermore,  the latter one is no bigger than
\[
\|f^{'}\|_{\infty,(r,R)}^{p}{\rm sup}\{ h(\lambda): \lambda \in [\epsilon/2,\|f\|_{\infty,(r,R) }+\epsilon/2]  \}.
\]
Therefore by the Lebesgue's Monotone and Dominated Convergence Theorems we have \[
A_{\eta} (r,R)\stackrel{\eta\rightarrow 0}{\longrightarrow}\int_{ (r,R)\cap \{x:f(x)>0\}}|f^{'}|^ph(f)dx.\]
 On the other hand,  we recognize that
\[
B_{\eta}(r,R) = \left(\sqrt{p-1}\right)^p\int_{E_\epsilon}\left|f^{''}{\cal T}_h(f_\eta)\right|^{\frac{p}{2}}h(f_\eta)dx+\left(\sqrt{p-1}\right)^p\int_{F_\epsilon}\left|f^{''}{\cal T}_h(f_\eta)\right|^\frac{p}{2}h(f_\eta)dx,
\]
because $f^{''}=0$ almost everywhere on the set $\{ x\in (r,R): f(x)=0\}$.
As before, we note that the former integrand increases as $\eta\rightarrow 0$, the latter one is no bigger than
 \[
 \left(\sqrt{p-1}\right)^p\| f^{''}\|_{\infty, (r,R)}^{p/2} {\rm sup}\{{|\cal T}_h(\lambda)|^{\frac{p}{2}}h(\lambda): \lambda\in [\epsilon/2,\|f\|_{\infty,(r,R) }+\epsilon/2]  \} .
 \]
 Therefore \[
 B_\eta (r,R)\stackrel{\eta\rightarrow 0}{\longrightarrow}\left(\sqrt{p-1}\right)^p\int_{(r,R)\cap \{x: f(x)>0\} }\left|f^{''}{\cal T}_h(f)\right|^\frac{p}{2}h(f)dx.\]
Now we verify the convergence of $\theta_\eta (r,R)$ as $\eta\to 0$. Note that if $f(R)\neq 0$ and $f(r)\neq 0$, we get $\theta_\eta (r,R)\rightarrow \theta_0  (r,R)$
as $\eta\to 0$.
In case $f(R)=0$ (respectively $f(r)=0$) we have $f^{'}(R)=0$ (respectively $f^{'}(r)=0$). Therefore in all cases $\theta_\eta (r,R)$ converges as $\eta\to 0$
to \[
\tilde{\theta}(r,R):= |f^{'}(R)|^{p-2}f^{'}(R) \tilde{H}(f(R)) - |f^{'}(r)|^{p-2}f^{'}(r) \tilde{H}(f(r)).
\]
 Altogether gives inequality (\ref{nier1111}) with the assumption $f\ge \eta$ relaxed to $f\geq 0$ almost everywhere and $\tilde{\theta} (r,R)$ instead of $\theta (r,R)$.
  To finish the proof of Case A, it suffices to chose the suitable sequence
 $r_n\to a$ and $R_n\to b$ applied for $r$ and $R$ and let $n$ converge to infinity.

 \smallskip
 \noindent
 {\sc Proof of Case B.}
Assume for simplicity that $h$ is bounded on some neighborhood of zero, while $|{\cal T}_h|^\frac{p}{2}h$ is nonincreasing on some neighborhood of zero.
The remaining cases follow by obvious modifications of our previous arguments. \\
We take  $\epsilon >0$ such that $h$ is bounded and  $|{\cal T}_h|^\frac{p}{2}h$ is nonincreasing on $(0,\epsilon]$ and repeat previous arguments up to inequality
(\ref{tez1}). Then we apply Lebegue's Dominated Convergence Theorem to both integrands on the right hand side of (\ref{tez1}) to deduce that $A_\eta (r,R)\stackrel{\eta\rightarrow 0}{\longrightarrow}\int_{(r,R)\cap \{ f>0\}}|f^{'}|^ph(f)dx$. The remaining arguments are the same as in Case A.
\noindent
This finishes the proof.\hfill$\Box$
\smallskip

\begin{rem}\label{roznef}\rm
If the assumptions of Proposition \ref{gorac} are satisfied and $f\in W^{2,1}_{loc}(a,b)$  is such that $f\ge 0$ on $(a,b)$, then we have
{\small
\begin{eqnarray*}
 \int_{(a,b)\cap \{ x: f(x)> 0\}}|f^{'}(x)|^ph(f(x))dx\le
  \left( \sqrt{p-1}\right)^p\int_{(a,b)\cap \{ x: f(x)> 0\}}\left|f^{''}(x){\cal T}_h(f(x))\right|^\frac{p}{2}h(f(x))dx\\ + \tilde{L}_{-}(a,b)(f),\ {\rm where}\\
\tilde{L}_{-}(a,b)(f) := \left( \liminf_{R\nearrow b}|f^{'}(R)|^{p-2}f^{'}(R)\tilde{H}(f(R)) - \limsup_{r\searrow a}|f^{'}(r)|^{p-2}f^{'}(r)\tilde{H}(f(r))\right).
\end{eqnarray*}
}
This follows from arguments in the proof of our last proposition.
\end{rem}

\section{Relaxing the nonnegativity assumption. Inequalities within less regular weights}\label{five}

 Our next goal is to examine when similar type of results hold with function $f$ not necessarily being nonnegative.
 One cannot apply directly Propositions \ref{goracs} and \ref{gorac} just to $|f|$ instead of $f$ as in such a case $|f|$ may not have the locally integrable distributional second derivative.

\noindent
 We start with the following easy observation stating that all results of our previous Chapter can be adapted to the case when function $f$ does not have
 single zeroes.

 \begin{prop}\label{doublezeroes}
 Assume that the following conditions hold:
 \begin{description}
 \item[1.] The assumptions on $p,a,b,h,H$ and ${\cal T}_h$ are as in Proposition \ref{goracs} and function $f$ belongs to the set
 {\small
\begin{eqnarray*}
\tilde{\cal R}^0_{dz}(a,b)&:=& \{ f\in W_{loc}^{2,1}((a,b)), f\ \hbox{\rm does not have single zeroes} : \\
&&  \liminf_{R\nearrow b}
|f^{'}(R)|^{p-2}f^{'}(R){\rm sign}(f(R))\tilde{H}(|f(R)|)\le 0,\ {\rm and}\\ &&
   \limsup_{r\searrow a}|f^{'}(r)|^{p-2}f^{'}(r){\rm sign}(f(r))\tilde{H}(|f(r)|) \ge 0
\},
\end{eqnarray*}
}
\item[2.] The assumptions on $p,a,b,h,H$ and ${\cal T}_h$ are as in Proposition \ref{gorac} and function $f$ belongs to the set
{\small
\begin{eqnarray*}
\tilde{\cal R}_{dz}(a,b)&:=& \{ f\in W_{loc}^{2,1}((a,b)),f\ \hbox{\rm does not have single zeroes} : \\
&&\left( \liminf_{R\nearrow b}|f^{'}(R)|^{p-2}f^{'}(R){\rm sign}(f(R))\tilde{H}(|f(R)|)\right.\\&&\left.
- \limsup_{r\searrow a}|f^{'}(r)|^{p-2}f^{'}(r){\rm sign}(f(r))\tilde{H}(|f(r)|)\right)\le 0
\},
\end{eqnarray*}
}
\end{description}
Then we have
{\small
\[
\int_{(a,b)\cap \{x: f(x)\neq 0\} }|f^{'}(x)|^ph(|f(x)|)dx\leq \left( \sqrt{p-1}\right)^p\int_{(a,b)\cap \{x: f(x)\neq 0\} }\left|f^{''}(x){\cal T}_h(|f(x)|)\right|^\frac{p}{2}h(|f(x)|)dx.
\]}
 \end{prop}

{\bf Proof.} This follows from the fact that if $f\in W^{2,1}_{loc}(a,b)$ does not have single zeroes then function $g=|f|$ belongs to $W^{2,1}_{loc}(a,b)$ and
$g^{'}=f^{'}{\rm sign}f$, $g^{''}=f^{''}{\rm sign}f$. Therefore  Propositions \ref{goracs} and \ref{gorac} can be applied to $|f|$ instead of $f$.\hfill$\Box$

Our next proposition applies to the situation when function $h$ possesses some regularity properties, namely when it is integrable in some (then also every) neighborhood of zero.

\begin{prop}\label{goracs2}
Let $p\geq 2$, $-\infty \le a<b\le+\infty$, $h:(0,\infty)\rightarrow(0,\infty)$ be a continuous function which is integrable in every neighborhood of zero, in particular
its primitive ${H}(x):=\int_0^x h(\tau)d\tau$ is well defined, extends  to $0$ and its extension $\tilde{H}$ satisfies $\tilde{H}(0)=0$.
 Let ${\cal T}_h:(0,\infty)\rightarrow(0,\infty)$ be as in Definition \ref{def1}.
Then for every function $f\in \tilde{\cal R}^0(a,b)$ where
{\small
\begin{eqnarray*}
\tilde{\cal R}^0(a,b)&:=& \{ f\in W_{loc}^{2,1}((a,b)):  \liminf_{R\nearrow b}|f^{'}(R)|^{p-2}f^{'}(R){\rm sign}f(R)\tilde{H}(|f(R)|)\le 0\\
&&  \ {\rm and}\  \limsup_{r\searrow a}|f^{'}(r)|^{p-2}f^{'}(r){\rm sign}f(r)\tilde{H}(|f(r)|) \ge 0
\},
\end{eqnarray*}
}
we have
{\small
\[
\int_{(a,b)\cap\{x:f(x)\neq 0\}}|f^{'}(x)|^ph(|f(x)|)dx\leq \left( \sqrt{p-1}\right)^p\int_{(a,b)\cap\{x:f(x)\neq 0\}}\left|f^{''}(x){\cal T}_h(|f(x)|)\right|^\frac{p}{2}h(|f(x)|)dx.
\]
}
\end{prop}

{\bf Proof.} We repeat the arguments of the proof of Proposition \ref{goracs} with $|f|=f{\rm sign}f$ instead of $f$ (so that ${\cal S}:= \{ x\in (a,b): |f(x)|>0\}$) and note that (the modified) formulas (\ref{sobota3}), (\ref{sobota1}), (\ref{sobota2})
hold as well because of an assumption $\tilde{H}(0)=0$. Then we complete the proof by the same arguments.  \hfill$\Box$

 Our last generalization strongly relays on the monotonicity properties of the involved functions. It reads as follows.

 \smallskip

 \begin{prop}\label{goal111}
Let $p\geq 2$ and $h:(0,\infty) \rightarrow (0,\infty)$ be a continuous function, $H, \ct_h:(0,\infty)\rightarrow\r$ be the same as in Definition \ref{def1}. Denote: $G_h(\lambda):=\frac{|\ct_h(\lambda)|^{p/2}h(\lambda)}{\lambda^{p/2}}$, ${\cal A}_f:= \{ r\in (a,b) : f^{'}(r)=0\ {\rm or}\ (f(r)\neq 0\  {\rm and}\ f^{'}(r)\neq 0)\}$ - the set of full measure\\ in $(a,b)$.  Assume further that:\\
{\bf (i)} both functions $h$ and $|\ct_h|^\frac{p}{2}h$ are either bounded or nonincreasing in some neighborhood of $0$;\\
{\bf (ii)} if $|\ct_h|^\frac{p}{2}h$ is not nonincreasing near zero then the function $G_h(\lambda)$ is either nonincreasing or bounded near $0$.
\noindent
Then for every $f\in {\cal R}_G(a,b)$, where
{\small
\begin{eqnarray*}
{\cal R}_G(a,b)&:=& \{ f\in W_{loc}^{2,1}(\r) : \exists_{\epsilon >0:}\forall_{r,R:a<r<R<b} \int_{ (r,R)\cap\{x: 0<|f(x)|<\epsilon\} }|f^{'}(x)|^pG_h(|f(x)|)dx <\infty  \\&&{\rm and}  \
\liminf_{R\nearrow b, R\in {\cal A}_f}|f^{'}(R)|^{p-2}f^{'}(R){\rm sign}f(R)H(|f(R)|)\chi_{\{ R:f(R)\neq 0\} }\\&& - \limsup_{r\searrow a, r\in {\cal A}_f}|f^{'}(r)|^{p-2}f^{'}(r){\rm sign}f(r)H(|f(r)|)\chi_{\{r: f(r)\neq 0\} } \le 0\}
\end{eqnarray*}
}
we have
{\small
\begin{equation*}
\int_{(a,b)\cap \{x:f(x)\neq 0\} }|f^{'}(x)|^ph(|f(x)|)dx\leq \left( \sqrt{p-1}\right)^p\int_{(a,b)\cap \{x:f(x)\neq 0\} }\left|f^{''}(x){\cal T}_h(|f(x)|)\right|^\frac{p}{2}h(|f(x)|)dx.
\end{equation*}
}
\end{prop}
\noindent
{\bf Proof.} The proof has the similar structure as that of Proposition \ref{gorac}.
\noindent
If\\ $\int_{ \{x: f^{'}(x)\neq 0\}}\left|f^{''}{\cal T}_h(|f|)\right|^\frac{p}{2}h(|f|)dx=\infty$ then obviously inequality holds.
Therefore we  suppose that this integral is finite and consider two cases:\\
\smallskip
{\sc Case A:} Both functions $h$ and $|\ct|^\frac{p}{2}h$ are nonincreasing in some neighborhood of zero;\\
{\sc Case B:} One of the functions $h$ or $|\ct_h|^\frac{p}{2}h$ is bounded in some neighborhood of zero. \\

{\sc Proof of Case A.}
We find some $\epsilon>0$ such that $h$ and $|\ct|^\frac{p}{2}h$ are nonincreasing on $(0,\epsilon]$ and denote $f_\eta (x):=\sqrt{f^2(x)+\eta}$ where $0<\eta<\frac{\epsilon}{2}$. An easy computation shows that
\begin{eqnarray*}
f_\eta^{'}(x) &=& \frac{f(x)f^{'}(x)}{f_\eta(x)},\\
f_\eta^{'}(x) &=&0 \ \ {\rm when}\ f^{'}(x)=0,\\
f_\eta^{''}(x) &=&\frac{f(x)}{f_\eta(x)}f^{''}(x) + \frac{(f^{'}(x))^2}{f_\eta(x)}\cdot \frac{\eta}{f^2(x)+\eta} \stackrel{\eta\to 0}{\rightarrow}{\rm sign}(f(x))f^{''}(x), \
{\rm when}\ f(x)\neq 0.
\end{eqnarray*}
According to Lemma \ref{pierwszylemat}
we may substitute $f_\eta$ to inequality (\ref{nier1111}), getting for every $a<r<R<b$
\begin{eqnarray*}
A_\eta(r,R)&:=&\int_r^R|f_\eta^{'}|^ph(f_\eta)dx\leq \left(\sqrt{p-1}\right)^p\int_r^R\left|f_\eta^{''}{\cal T}_h(f_\eta)\right|^\frac{p}{2}h(f_\eta)dx+\theta_\eta(r,R)\label{pomoc}\\
&=:&B_\eta(r,R)+\theta_\eta(r,R),
\end{eqnarray*}
where \[
\theta_\eta (r,R) := |f_\eta^{'}(R)|^{p-2}f_\eta^{'}(R) H(f_\eta(R)) - |f_\eta^{'}(r)|^{p-2}f_\eta^{'}(r) H(f_\eta(r)).\]
We note that $f_\eta\searrow|f|$ and $|f_\eta^{'}|\nearrow|f^{'}|$ when $\eta\searrow 0$. Therefore by almost the same argument as in the proof of Proposition \ref{gorac},  we have \[
A_{\eta} (r,R)\stackrel{\eta\rightarrow 0}{\longrightarrow}\int_{\{x\in (r,R):  f(x)\neq 0\}}|f^{'}|^ph(|f|)dx.\]
On the other hand, we have
\begin{eqnarray*}
B_{\eta}(r,R)=\left(\sqrt{p-1}\right)^p\int_{\{x\in (r,R): f(x)\neq 0\}}\left|f_\eta^{''}{\cal T}_h(f_\eta)\right|^{\frac{p}{2}}h(f_\eta)dx .
\end{eqnarray*}
 Moreover $|f^{''}_\eta\ct_h(f_\eta)|^{p/2}h(f_\eta)\stackrel{\eta\rightarrow 0}{\longrightarrow}|f^{''}\ct_h(|f|)|^{p/2}h(|f|)$ a.e. on $\{x\in(r,R): f(x)\neq 0\}$. In order to obtain convergence of integrals we will estimate the integrand of $B_\eta(r,R)$ and apply  Lebesgue's Dominated Convergence Theorem. We note that
\[
|f_\eta^{''}|\le |f^{''}|+\frac{|f^{'}|^2}{f_\eta},\  {\rm  on}\  [r,R].
\]
Therefore
\[
\left|f^{''}_\eta\ct_h(f_\eta)\right|^\frac{p}{2}h(f_\eta)\leq I_\eta +II_\eta,
\]
where
\begin{eqnarray*}
I_\eta := \left|2f^{''}{\cal T}_h(f_\eta )\right|^{\frac{p}{2}}h(f_\eta ), &
II_\eta  := \left|\frac{2|f^{'}|^2}{f_\eta}{\cal T}_h(f_\eta)\right|^{\frac{p}{2}}h(f_\eta).
\end{eqnarray*}
Now we estimate $I_\eta$ and $II_\eta$. We have for almost every $x\in [r,R]$:
\begin{eqnarray*}
I_\eta=\left|2f^{''}{\cal T}_h(f_\eta)\right|^{\frac{p}{2}}h(f_\eta)\chi_{E_\epsilon} + \left|2f^{''}{\cal T}_h(f_\eta)\right|^{\frac{p}{2}}h(f_\eta)\chi_{F_\epsilon},
\end{eqnarray*}
where $E_\epsilon=\{x\in (r,R):   0<|f(x)|<\frac{\epsilon}{2}\}$ and $F_\epsilon=\{x\in (r,R):  |f(x)|\ge \frac{\epsilon}{2}\}$. The former function increases as $\eta\rightarrow 0$ and is no bigger than
$2^\frac{p}{2}\left|f^{''}\ct_h(|f|)\right|^\frac{p}{2}h(|f|)$, which is an integrable function. The latter function is no bigger than
\[
 2^\frac{p}{2}\| f^{''}\|_{\infty, (r,R)}^{p/2} {\rm sup}\{|{\cal T}_h(\lambda)|^{\frac{p}{2}}h(\lambda): \lambda\in [\epsilon/2,\|f\|_{\infty,(r,R) }+\epsilon/2]  \} ,
 \]
which is a  constant, so that it is also integrable on $[r,R]$. On the other hand, we have for almost every $x\in [r,R]$:
\begin{eqnarray*}
II_\eta=\left|\frac{2|f^{'}|^2}{f_\eta}{\cal T}_h(f_\eta)\right|^{\frac{p}{2}}h(f_\eta)\chi_{E_\epsilon}+ \left|\frac{2|f^{'}|^2}{f_\eta}{\cal T}_h(f_\eta)\right|^{\frac{p}{2}}h(f_\eta)\chi_{F_\epsilon}.
\end{eqnarray*}
The former function increases as $\eta\rightarrow 0$ to $2^{p/2}|f^{'}|^{p}G_h(|f|)\chi_{\{x:0<|f(x)|<\epsilon\}}$, which by assumption is an integrable function. The latter function is no bigger than
\begin{eqnarray*}
2^\frac{p}{2}\|f^{'}\|^p_{\infty,(r,R)}{\rm sup}\{G_h(\lambda): \lambda\in[\epsilon/2, \|f\|_{\infty,(r,R)}+\epsilon/2]<\infty\}.
\end{eqnarray*}   Thus by the Lebesgue's Dominated Convergence Theorem we obtain
\begin{eqnarray*}
B_{\eta} (r,R)\stackrel{\eta\rightarrow 0}{\longrightarrow}\left( \sqrt{p-1}\right)^p\int_{\{x\in(r,R):f(x)\neq 0\}}\left|f^{''}(x){\cal T}_h(|f(x)|)\right|^\frac{p}{2}h(|f(x)|)dx.
\end{eqnarray*}
Now we verify the convergence of $\theta_\eta(r,R)$ as $\eta\rightarrow 0$. By almost the same arguments as in the proof of Proposition \ref{gorac},  the expression $|f_\eta^{'}(R)|^{p-2}f_\eta^{'}(R)H(f_\eta(R))$ might not converge to $|f^{'}(R)|^{p-2}f^{'}(R)H(|f(R)|)$ only if
$R\in \{(a,b): f(R)=0, f^{'}(R)\neq 0\} = (a,b)\setminus {\cal A}_f$, which is the set of measure zero in $(a,b)$.
Therefore for  every $R\in {\cal A}_f$,  we have
\[
 \lim_{\eta\to 0} |f_\eta^{'}(R)|^{p-2}f_\eta^{'}(R)H(f_\eta (R))\rightarrow |f^{'}(R)|^{p-2}f^{'}(R){\rm sign}f(R)H(|f (R)|)\chi_{\{f(R)\neq 0\}}.
 \]
By the same arguments, for almost every $r\in {\cal A}_f$
\[
 \lim_{\eta\to 0} |f_\eta^{'}(r)|^{p-2}f_\eta^{'}(r)H(f_\eta (r))\rightarrow |f^{'}(r)|^{p-2}f^{'}(r){\rm sign}f(r)H(|f (r)|)\chi_{\{f(r)\neq 0\}}.
 \]
Altogether gives that for every $r,R\in {\cal A}_f$:
 \begin{eqnarray*}
  \int_r^R|f^{'}(x)|^ph(|f(x)|)dx&\leq& \left( \sqrt{p-1}\right)^p\int_r^R\left|f^{''}(x){\cal T}_h(|f(x)|)\right|^\frac{p}{2}h(|f(x)|)dx +\tilde{\theta} (r,R),
 \end{eqnarray*}
 where
\begin{eqnarray*}
 \tilde{\theta}(r,R) =& |f^{'}(R)|^{p-2}f^{'}(R){\rm sign}f(R) H(|f(R)|)\chi_{\{R:f(R)\neq 0\}}-\\ & |f^{'}(r)|^{p-2}f^{'}(r){\rm sign}f(r) H(|f(r)|)\chi_{\{r:f(r)\neq 0\}}.
\end{eqnarray*}
It remains to chose a suitable sequence $r_n\rightarrow a$ and $R_n\rightarrow b$, $r_n, R_n\in{\cal A}$ and let $n$ converge to infinity. This finishes the proof of Case A.\\
\smallskip
{\sc Proof of Case B.}  We proceed similarly to the proofs of Proposition \ref{gorac} and Case A. The only noteworthy difference appears when $|\ct_h|^\frac{p}{2}h$ is bounded (and is not decreasing) near zero since we cannot estimate $II_\eta$ analogously to any of the previous cases. However, in this case, $G_h$ is either nonincreasing or bounded near zero. For nonincreasing $G_h$ we estimate $II_\eta$ exactly like in Case A, while  for bounded $G_h$ we note that $II_\eta$ is no bigger than
\[
2^\frac{p}{2}\|f^{'}\|_{\infty,(r,R)}^{p}{\rm sup}\{G_h(\lambda):\lambda\in(0,\|f\|_{\infty,(r,R)})\},
\]
which is a positive constant. This finishes the proof.\hfill$\Box$

\begin{rem}\rm
Assumption $\int_{ (r,R)\cap\{ 0<|f|<\epsilon\} } G_g(|f(x)|)dx <\infty$ for all $a<r<R<b$ cannot be omitted in general from the definition of the class of admitted functions ${\cal R}_G(a,b)$. Counterexample will be provided in Remark \ref{trzykoty} in the next section.
\end{rem}

\section{Examples}\label{six}
In the sequel we present some examples illustrating Propositions \ref{goal2}, \ref{gorac}, \ref{doublezeroes}, \ref{goracs2} and \ref{goal111}.
For simplicity we assume that the admitted function $f\in W_{loc}^{2,1}(\r)$ has a compactly supported first weak derivative.
Obviously, this assumption can be essentially weakened, taking into account our more general statements.
The discussion is provided by four propositions stated below, dealing with power, logarithmic and exponential type functions $h$.

\begin{prop}\label{przyk1} Let $2\le p<\infty$, $\theta\in\r$ and $f\in W_{loc}^{2,1}(\r)$ such that $f^{'}$ has a compact support. Moreover, suppose
that at least one of the following assumptions is satisfied:
\begin{enumerate}
\item $\theta<\frac{1}{p}$,
\item $\theta >\frac{1}{p}$ and $f$ is nonnegative or (more generally) does not have single zeroes,
\item $ \theta > \frac{1}{p}$ and there exists $\epsilon$ such that for all  $r<R$: $ \int_{(r,R)\cap\{x:0<|f(x)|<\epsilon\}}\left(\frac{|f^{'}|}{|f|^\theta}\right)^pdx<\infty $.
\end{enumerate}
Then we have
\begin{eqnarray}\label{ex1}
\int_{\{x:f(x)\neq 0\}} \left(\frac{ |f^{'}|}{|f|^\theta}\right)^pdx \le \left( \frac{p-1}{|1-\theta p|}\right)^{\frac{p}{2}}
 \int_{\{x:f(x)\neq 0\}} \left(\frac{\sqrt{ |ff^{''}|}}{|f|^\theta}\right)^{{p}}dx.
\end{eqnarray}
\end{prop}
\noindent

{\bf Proof.} We set: $h(\lambda)= \frac{1}{\lambda^{\theta p}}$. It  is a continuous,  positive function defined on $(0,\infty)$ and we choose it's primitive: $H(\lambda)=\frac{{\lambda}^{1-\theta p}}{1-\theta p}$. Then $|{\cal T}_h(\lambda)|^\frac{p}{2}h(\lambda)=\frac{\lambda^{(\frac{p}{2}-\theta p)}}{|1-\theta p|^\frac{p}{2}}$. Depending on $\theta$,  the function $|{\cal T}_h(\lambda)|^\frac{p}{2}h(\lambda)$ is locally bounded ($\theta\leq\frac{1}{2}$) or nonincreasing ($\theta>\frac{1}{2}$).
\\
1:  In case  $\theta\leq 0$ let $\alpha=-\theta$, so that $h(\lambda)=\lambda^{\alpha p}:(0,\infty)\rightarrow(0,\infty)$. As then $H(\lambda)=\int_0^\lambda h(s)ds=\frac{\lambda^{\alpha p+1}}{\alpha p+1}$, it  extends continuously to $0$ and it's extension ${\tilde H}$ satisfies: ${\tilde H}(0)=0$. Therefore
 we can apply Proposition \ref{goracs2} to get
\begin{eqnarray*}
\int_\r\left(|f^{'}||f|^\alpha\right)^pdx\leq\left(\frac{p-1}{\alpha p+1}\right)^\frac{p}{2}\int_\r\left|f^{''}f\right|^\frac{p}{2}\cdot|f|^{\alpha p}dx.
\end{eqnarray*}\\
Almost the same arguments allow to apply Proposition \ref{goracs2} when $0\leq\theta<\frac{1}{p}$.

\noindent
2:  Assume  that $\theta >\frac{1}{p}$ and $f$ is nonnegative. Then $h$ is nonincreasing and   by Proposition \ref{gorac} we have
\begin{eqnarray*}
\int_{\{x:f(x)\neq 0\}} \left(\frac{ |f^{'}|}{f^\theta}\right)^pdx\leq\left(\sqrt{p-1}\right)^p\int_{\{x:f(x)\neq 0\}}\left(|f^{''}|\frac{f}{|1-\theta p|}\right)^\frac{p}{2}\cdot \frac{1}{f^{\theta p}} dx=\\
=\left(\frac{p-1}{|1-\theta p|}\right)^\frac{p}{2}\int_{\{x:f(x)\neq 0\}}\left(\frac{\sqrt{ |ff^{''}|}}{f^\theta}\right)^{{p}}dx.\\
\end{eqnarray*}
If $f$ does not have single zeroes the only difference  is that now we apply Proposition \ref{doublezeroes}, part 2.
\\
3: In the last case we  apply Proposition \ref{goal111}. Now function $h$ is nonincreasing, however we also know that for $\theta\leq\frac{1}{2}$ the function $|\ct_h|^\frac{p}{2}h$ is not nonincreasing. Therefore we require the function $G_h(\lambda)=\frac{|\ct_h(\lambda)|^\frac{p}{2}h(\lambda)}{\lambda^\frac{p}{2}}$ to be either nonincreasing or bounded in some neighborhood of zero.  In our case $G_h(\lambda)=\frac{1}{|1-\theta p|^\frac{p}{2}|\lambda|^{\theta p}}$  is nonincreasing. Thus the inequality holds for all $f$ satisfying our assumptions.
This finishes the proof.$\hfill\Box$

\begin{rem}\rm
The case with $\theta =0$ was already obtained in Lemma \ref{motyw}. The case with $\theta =\frac{1}{2}$ (then $f$ on the right hand side does not appear) and nonnegative $f\in C_0^\infty(\r)$ can be found in \cite{ma}, Lemma 1, Section 8.2.1.
\end{rem}

\noindent
Now we deal with case $p=\infty$ obtained as the limiting case from Proposition \ref{przyk1}.

\begin{prop}\label{przyk2}
Let  $\theta\in\r$ and $f\in W_{loc}^{2,1}(\r)$ such that $f^{'}$ has a compact support. Moreover, suppose
that at least one of the following assumptions is satisfied:
\begin{enumerate}
\item $\theta\le 0$,
\item $\theta >0$ and $f$ is nonnegative or  does not have single zeroes,
\item $ \theta >0$ and there exists $\epsilon$ such that for all  $r<R$: $ \int_{(r,R)\cap\{x:0<|f(x)|<\epsilon\}}\left(\frac{|f^{'}|}{|f|^\theta}\right)^pdx<\infty $, whenever $2\le p <\infty$.
\end{enumerate}
Then we have
\[
\left\|\frac{f^{'}}{|f|^\theta}\chi_{\{x:f(x)\neq 0\}}\right\|_\infty\leq\frac{1}{\sqrt{\theta}} \left\|\frac{\sqrt{|ff^{''}|}}{|f|^\theta}\chi_{\{x:f(x)\neq 0\}}\right\|_\infty.
\]
\end{prop}

\noindent
{\bf Proof.} In all considered cases we find finite $p_0$ such that $f$ satisfies either: 1, 2 or 3 in  Proposition \ref{przyk1}, with all $p\ge p_0$.
Consequently:
\begin{eqnarray*}
\left( \int_{\{x:f(x)\neq 0\}} \left(
\frac{ |f^{'}| }{|f|^\theta}
\right)^pdx\right)^{\frac{1}{p}} \le \left( \frac{p-1}{|1-\theta p|}\right)^{\frac{1}{2}}
 \left( \int_{\{x:f(x)\neq 0\}} \left(\frac{\sqrt{ |ff^{''}|}}{|f|^\theta}\right)^{{p}}dx\right)^{\frac{1}{p}}
\end{eqnarray*}
for every finite $p\ge p_0$. Now it suffices to let $p\to\infty$.\hfill$\Box$

\begin{rem}[necessity of the assumptions]\label{trzykoty}\rm
Let us discuss Assumption 2 on $f$ in Proposition \ref{przyk1}: {\bf (A2):}{\it $f$ is nonnegative or does not have single zeroes}.
To show necessity of this assumption we consider $f(x)=\phi(x)\sin 2\pi x$, where $\phi\in C_0^\infty(\r)$, $\phi\equiv 1$ on $[-\frac{1}{2},\frac{1}{2}]$, $\phi\equiv 0$ on $(-\infty,-1)\cup(1,\infty)$, $\phi^{'}\neq 0$ on $(-1,-1/2)\cup (1/2,1)$, $0\le\phi\le 1$. It is clear that $f$ changes its sign at $0$ and has single zero at $0$, so it does not satisfy this assumption. Moreover,
\begin{eqnarray*}
f^{'}(x)&=&\phi^{'}(x)\sin 2\pi x+\phi(x)2\pi\cos 2\pi x,\\
f^{''}(x)&=&\phi^{''}(x)\sin 2\pi x + 2\phi^{'}(x)2\pi\cos 2\pi x+ \phi(x)(2\pi)^2\sin 2\pi x,
\end{eqnarray*}
so that for almost every $x$
\[
\frac{|f^{'}(x)|}{|f(x)|^\theta}\chi_{\{f(x)\neq 0\}}\ge \frac{|f^{'}(x)|}{|f(x)|^\theta}\chi_{(-1/8,1/8)} =\left|\frac{2\pi \cos 2\pi x}{\sin 2\pi x}\right|\chi_{(-1/8,1/8)}\ge
\frac{\pi\sqrt{2}}{|x|^\theta}\chi_{(-1/8,1/8)}.
 \]
 Therefore when $\theta >\frac{1}{p}$ left hand side in (\ref{ex1}) is infinite. Furthermore, when $\theta \le \frac{1}{2}$ function $\frac{\sqrt{ |ff^{''}|}}{|f|^\theta}$ is bounded, so $p$-integrable (as compactly supported). Therefore right hand side in  (\ref{ex1}) is finite. It shows necessity of the assumption {\bf (A2)}
 in Proposition \ref{przyk1} in general.  The same example shows necessity of assumption {\bf (A3):} {\it there exists $\epsilon$ such that for all  $r<R$: $ \int_{(r,R)\cap\{x:0<|f(x)|<\epsilon\}}\left(\frac{|f^{'}|}{|f|^\theta}\right)^pdx<\infty $} in the considered range of parameters.
 Moreover, assumptions:  {\bf (A2)} and {\bf (A3)} appear also to be necessary in the limiting case: $p=\infty$ in Assumptions 2 and 3 of Proposition
 \ref{przyk2}. The argument is provided by the same example.
\end{rem}

\noindent
Our next proposition complements the case $\theta =\frac{1}{p}$ in the statement of Proposition \ref{przyk1}.

  \begin{prop}\label{przyk3}
  Let $2\le p <\infty$ and $f\in W_{loc}^{2,1}(\r)$ such that $f^{'}$ has a compact support. Moreover suppose $f$ is nonnegative or does not have single zeroes or there exists $\epsilon$ such that for all $r<R$ we have $\int_{(r,R)\cap\{x:0<|f(x)|<\epsilon\}}\frac{|f^{'}|^p|\ln|f||^{p/2}}{|f|}dx<\infty.$
Then we have
\begin{eqnarray*}
\int_{\{x:f(x)\neq 0\}}\frac{|f^{'}|^p}{|f|}dx\leq \left(\sqrt{p-1}\right)^p\int_{\{x:f(x)\neq 0\}}\frac{|ff^{''}\ln(|f|)|^\frac{p}{2}}{|f|}dx.
\end{eqnarray*}
\end{prop}

\noindent
{\bf Proof.}  Function $h(\lambda)=\frac{1}{\lambda}$ is  continuous, nonincreasing and positive on $(0,\infty)$.  We choose it's primitive $H(\lambda)=\ln\lambda$. Then the function $|\ct_h(\lambda)|^\frac{p}{2}h(\lambda)=\lambda^{\frac{p}{2}-1}|\ln\lambda|^\frac{p}{2}$ is locally bounded for $p>2$. On the other hand, for $p=2$,  we have $|\ct_h(\lambda)|h(\lambda)=|\ln\lambda|$, it is nonincresing near zero. Either way, when $f$ is nonnegative or does not have single zeroes, we can apply Proposition \ref{gorac} or its extension, Proposition \ref{doublezeroes}.
 If this is not the case we aim to apply Proposition \ref{goal111}. We already know that $h$ and $|\ct_h(\lambda)|^\frac{p}{2}h(\lambda)$ satisfy required  assumptions.
  Therefore it remains to show that if $p>2$ (in particular $|\ct_h(\lambda)|^\frac{p}{2}h(\lambda)$ is not nonincreasing near zero) then the function $G_h(\lambda)=\frac{|\ln\lambda|^{p/2}}{\lambda}$ is either bounded or nonincreasing near zero. Crearly it's nonincreasing near zero. Therefore the inequality holds for every $f$ satisfying the required assumptions.\hfill$\Box$

Our last proposition provides example dealing with exponential function $h$.

 \begin{prop}\label{przyk4}
Let $p\geq 2$ and $f\in W^{2,1}_{loc}(\r)$ such that $f^{'}$ has a compact support. Furthermore, suppose that $\alpha,\beta\in\r$, $\alpha\neq 0$, $\beta\geq 1$, and one of the following conditions is satisfied:
\begin{enumerate}
\item $f$
is nonnegative or does not have single zeroes or there exists $\epsilon>0$ such that for all $r<R$ we have,\\ $\int_{(r,R)\cap\{x:0<|f(x)|<\epsilon\}}|f^{'}|^p|f|^{-\beta(\frac{p}{2}-1)-1}e^{\alpha|f|^\beta}dx<\infty$.
\item $\alpha> 0$ or $\beta=1$.
\end{enumerate}
Then we have
\begin{eqnarray}\label{czterykoty}
\int_{\{x:f(x)\neq 0\}}|f^{'}|^p|f|^{\beta-1}e^{\alpha |f|^\beta}dx\leq\left(\sqrt{\frac{p-1}{|\alpha|\beta}}\right)^p \int_{\{x:f(x)\neq 0\}}|f^{''}|^\frac{p}{2}|f|^{(\frac{p}{2}-1)(1-\beta)}e^{\alpha |f|^\beta}dx.
\end{eqnarray}
\end{prop}

\noindent
{\bf Proof.} 1. Suppose that $f$ is nonnegative.  We substitute in Proposition \ref{gorac}:  $h(\lambda):=\lambda^{\beta-1} e^{\alpha\lambda^\beta}$. Then the function $h:(0,\infty)\rightarrow(0,\infty)$ is continuous and bounded near $0$. We choose $H(\lambda)=\frac{1}{\alpha\beta}e^{\alpha\lambda^\beta}$ and notice that the function $|\ct_h(\lambda)|^\frac{p}{2}h(\lambda)=\frac{1}{(\sqrt{|\alpha|\beta})^p} \lambda^{(\frac{p}{2}-1)(1-\beta)}e^{\alpha\lambda^\beta}$ is positive and nonincreasing in a neighborhood of zero. For $f$ without single zeroes we apply Proposition \ref{doublezeroes}, part 2. In the last case we apply Proposition \ref{goal111}. We repeat previous arguments and observe that since $h$ is bounded and $|\ct_h|^\frac{p}{2}h$ is nonincreasing near $0$, the assumptions of Proposition \ref{goal111} are satisfied. Then we compute  that $G(\lambda)=\left(\frac{1}{|\alpha|\beta}\right)^\frac{p}{2}\lambda^{-\beta(\frac{p}{2}-1)-1} e^{\alpha\lambda^\beta}$.
 Therefore   inequality (\ref{czterykoty})  holds for all $f$ satisfying our assumptions.\\
2. Suppose that $\alpha>0$. We proceed similarly to the case 1. and observe that $H_1(\lambda):= H(\lambda)-\frac{1}{\alpha\beta}=\frac{1}{\alpha\beta}(e^{\alpha\lambda^\beta}-1)$ is an absolutely continuous primitive of $h$ that can be continuously extended to $0$ and it's extension ${\tilde H}_1$ satisfies:
${\tilde H}_1(0)=0$. Therefore we can apply Proposition \ref{goracs2} to obtain
\begin{eqnarray*}
\int_{\{x:f(x)\neq 0\}}|f^{'}|^ph(|f|)dx\leq\int_{\{x:f(x)\neq 0\}}\left|f^{''}\frac{H(|f|)-c}{h(|f|)}\right|^\frac{p}{2}h(|f|)dx,
\end{eqnarray*}
where $c=\frac{1}{\alpha\beta}$. To finish the proof it suffices to notice that \[
|H_1(|f|)|=|H(|f|)-c|=H(|f|)-c<H(|f|).\]
If $\beta=1$  we substitute $h(\lambda):=e^{\alpha\lambda}$, which is the continuous, positive function defined on $\r$ and apply  Proposition \ref{goal2}. $\hfill\Box$

\section{Applications to nonlinear eigenvalue problems}\label{seven}

\subsection{Explanation of general approach}

{\bf Eigenvalue problem.}
Suppose that $-\infty\le a<b\le +\infty$, $\tau: (A,B)\rightarrow \r$ is continuous, where $-\infty\le A<B\le\infty$, $g\in L^q(a,b)$ with $q\in [1,\infty ]$ and
 $f\in  W^{2,1}_{loc}((a,b))$, $f(x)\in (A,B)$ a.e. and  satisfies the following ODE:
\begin{eqnarray}\label{eq}
\left\{
\begin{array}{ccc}
f^{''}(x)&=&g(x)\tau(f(x))\ \hbox{\rm a.e. in}\ (a,b),\\
f\in {\cal R}
\end{array}
\right.
\end{eqnarray}
where set ${\cal R}$ will serve to define boundary conditions and will be indicated later.
Our goal is to establish regularity results for solutions to (\ref{eq}).

\noindent
{\bf Our approach.} For simplicity we only analyze the case when $f$ is nonnegative and set
 set $\{ \lambda\in (A,B): \tau(\lambda)=0\}$ is at most countable. Then equation (\ref{eq}) is automatically satisfied almost everywhere on the set where $\{ x:\tau (f(x))=0\}$
 as on level sets of $f$ we have $f^{''}=0$ a.e..
By  formal computation we find  such function $h$  that
\begin{equation}\label{guzik}
|g(x)|^q=\left|\frac{f^{''}(x)}{\tau(f(x))}\right|^q=|\ct_h(f(x))f^{''}(x)|^\frac{p}{2}h(f(x)),\ {\rm when}\ \tau(f(x))\neq 0,
\end{equation}
where $p=2q$ and $\ct_h$ is as in Definition \ref{def1}. This can be obtained by looking for such $h\geq 0$  that
\begin{equation}\label{oby}
|\tau(\lambda)|^{-q}=|\ct_h (\lambda)|^qh(\lambda)=|H(\lambda)|^qh^{1-q}(\lambda), \ \hbox{\rm provided that}\ \tau (\lambda)\neq 0.
\end{equation}
For this, let us consider two cases: $q>1$ and $q=1$ separately.

When $q>1$  function $k(\lambda):=|\tau(\lambda)|^\frac{q}{q-1}$ is well defined and locally $L^1$ on $(A,B)$.
Therefore it possesses the  locally absolutely continuous  primitive denoted by $K$, so that $K^{'}=k$.
Let us additionally assume that $K$ can be chosen in such a way that $K(\lambda)\neq 0$ when $\tau (\lambda )\neq 0$, i. e. $K$ has constant sign
on each component of set $\{ \lambda\in (A,B): \tau(\lambda)\neq 0\}$.
 An easy verification
shows that \[
H(\lambda):=
 -{\rm sgn K(\lambda )}(q-1)^{q-1}|K (\lambda)|^{1-q}
\]
satisfies $H^{'}=h$,  where
$h(\lambda):=(q-1)^{q}|K (\lambda)|^{-q}k.$ Moreover,
function $h$ is well defined when $\tau (\lambda)\neq 0$ and equations (\ref{oby}) and (\ref{guzik}) are satisfied.\\
When $q=1$ let us assume additionally that function $\tau$ has a continuous derivative on set  $\{ \lambda\in (A,B): \tau(\lambda)\neq 0\}$ and that $\tau^{'}$ has
constant sign on each component of $\{ \lambda\in (A,B): \tau(\lambda)\neq 0\}$. Then function
\[
H(\lambda) =-\frac{{\rm sign}\tau^{'}(\lambda)}{\tau(\lambda)}
\]
satisfies $H^{'}=h :=\frac{|\tau^{'}|}{\tau^2}$ on set $\{ \lambda\in (A,B): \tau(\lambda)\neq 0\}$. Therefore  equations (\ref{oby}) and (\ref{guzik}) are also satisfied.

Suppose further that $h$ and $f$ obey the assumptions of one of the Propositions: \ref{goal2}, \ref{winosekdodatk}, \ref{goracs}, \ref{gorac} respectively. Then we can deduce that
\[
\int_{(a,b)\cap\{x:f(x)\neq 0\}}|f^{'}|^{2q}h(f)\leq C\|g\|_q^q,
\]
in particular $|f^{'}|^{2q}h(f)\in L^1(a,b)$ and we can deduce further regularity results on $f$. Obvious modifications
and usage of Propositions: \ref{doublezeroes}, \ref{goracs2} or \ref{goal111} respectively allows to deal with nonlinearities like $\tau (|f|)$ instead
of $\tau (f)$ and obtain similar type of results.

 \subsection{ODEs dealing with strictly positive functions}

This subsection is devoted to strictly positive solutions of (\ref{eq}), the case when $(A,B)=(0,\infty)$. It is well known that such restriction
arises naturally in many nonlinear phenomena describing physical models (\cite{ar,arw}), as well as for example in the pure mathematical approach in the study of Poincare-type inequalities in functional analysis (see e.g. \cite{lind}).

\subsubsection{General result}

For this purpose we define the following three transforms of $\tau$.

\begin{defi}\label{def2}
Let $q\geq 1$, $\tau:(0,+\infty)\rightarrow(0,+\infty)$ be a continuous function. Furthermore assume that:\\
\smallskip
\noindent
a) if $q=1$ then $\tau$ has a continuous derivative which is of constant sign,\\
b) if $q>1$ function  $k_{q,\tau}(\lambda)=(\tau(\lambda))^\frac{q}{q-1}:(0,+\infty)\rightarrow(0,+\infty)$   admits a locally absolutely continuous  primitive denoted by $K_{q,\tau}$ which is of constant sign. \\
\smallskip
\noindent
We define the following three transforms of $\tau$
\begin{eqnarray*}
h_{q,\tau}(\lambda)&=&\left\{
\begin{array}{ccc}
(q-1)^q |K_{q,\tau}(\lambda)|^{-q}k_{q,\tau}(\lambda) &{\rm if}& q>1\\
\frac{|\tau^{'}(\lambda)|}{\tau(\lambda)^2}& {\rm if} &q=1
\end{array}
\right. \\
H_{q,\tau}(\lambda)&=&\left\{
\begin{array}{ccc}
-({\rm sgn}\, K_{q,\tau}(\lambda))(q-1)^{q-1}|K_{q,\tau}(\lambda)|^{1-q} &{\rm if}& q>1\\
-\frac{{\rm sgn \tau^{'}(\lambda) }}{\tau (\lambda)} & {\rm if} &q=1
\end{array}
\right. \\
G_{q,\tau}(\lambda)&=&\int h_{q,\tau}(\lambda)^\frac{1}{2q}d\lambda
\end{eqnarray*}
where $\int h_{q,\tau}(\lambda)^\frac{1}{2q}d\lambda$ is a locally absolutely continuous primitive of $(h_{q,\tau})^\frac{1}{2q}$.
\end{defi}

Applying our arguments from the beginning of this section we immediately obtain the following fact.

\begin{lem}\label{lemat} If $f\in W^{2,1}_{loc}(a,b)$ is a positive solution of (\ref{eq}) and $h=h_{q,\tau}$ is as in Definition \ref{def2},   then
equation
\begin{equation}\label{srod}
|g(x)|^q=\left|\frac{f^{''}(x)}{\tau(f(x))}\right|^q=|\ct_h(f(x))f^{''}(x)|^{q}h(f(x)),
\end{equation}
 is satisfied for almost every $x\in (a,b)$, where $\ct_h(f(x))$ is as in Definition \ref{def1} (with $A=0, B=\infty$).
\end{lem}

We arrive at the following result.

\begin{prop}\label{regul0}
Suppose that $1\le q<\infty$, $-\infty\leq a<b\leq +\infty$, $g\in L^q(a,b)$, $\tau:(0,+\infty)\rightarrow(0,+\infty)$ is continuous.
Further, let $h=h_{q,\tau}, H=H_{q,\tau}, G=G_{q,\tau}$ be as in Definition \ref{def2}
and let $f\in W^{2,1}_{loc}(a,b)$ be a positive solution of the following ODE:
\begin{equation*}
\left\{
\begin{array}{c}
f^{''}(x)=g(x)\tau(f(x))\  \hbox{\rm a.e. on}\ (a,b)\\
 \liminf_{R\nearrow b}|f^{'}(R)|^{2q-2}f^{'}(R){H}(f(R))
 -  \limsup_{r\searrow a}|f^{'}(r)|^{2q-2}f^{'}(r){H}(f(r)) \le 0
 \end{array}
 \right.
\end{equation*}
Then we have
\begin{description}
\item[i)] $|f^{'}|^{2q}h(f)\in L^1(a,b)$ and
\begin{eqnarray*}
\int_a^b|f^{'}(x)|^{2q}h(f(x))dx\le \left( \sqrt{2q-1}\right)^{2q}\int_a^b |g(x)|^qdx.
\end{eqnarray*}
\item[ii)] $F(x) := G(f(x))$ is  Holder continuous with exponent $\gamma=1-\frac{1}{2q}$ and
\begin{eqnarray*}
%\frac{|F(x)-F(y)|}{|x-y|^{1-\frac{1}{2q}}}\le \left( \int_x^y |f^{'}(x)|^{2q}h(f(x))dx \right)^{\frac{1}{2q}},\ \hbox{\rm when}\ a<x<y <b\nonumber\\
{\rm sup}\left\{  \frac{|F(x)-F(y)|  }{|x-y|^{1-\frac{1}{2q}}} : {x,y\in (a,b)}\right\} \le  \sqrt{2q-1} \left( \int_a^b |g(x)|^qdx\right)^{\frac{1}{2q}}.
\end{eqnarray*}
\item[iii)]
If $G$ is of constant sign and $|G|$ is
 strictly increasing then $f$ satisfies the following estimation:
\begin{eqnarray*}
|f(x)|\le |G|^{-1}\left( |G({f}(c)| +|x-c|^{1-\frac{1}{2q}}\sqrt{2q-1}\left( \int_a^b |g(\tau)|^qd\tau)\right)^{\frac{1}{q}}     \right), \\
\end{eqnarray*}
where $x,c\in (a,b)$ can be chosen arbitrary.
In particular when $a,b\in\r$, we have
\begin{eqnarray*}
|f(x)| \le G^{-1}\left( G({f}(c)) +({\max} \{ |c-a|, |b-c|\})^{1-\frac{1}{2q}}\sqrt{2q-1}\left( \int_a^b |g(\tau)|^qd\tau)\right)^{\frac{1}{q}}     \right),
\end{eqnarray*}
whenever  $x,c\in (a,b)$.
\item[iv)] If $a\in \r$, $G$ is of constant sign, $|G|$ is
 strictly increasing and $\lim_{r\to\infty}|G(r)|=\infty$,  then
$f$ extends to a continuous function up to $a$ (denote this extension by $\tilde{f}$) and we have
\begin{eqnarray*}
|f(x)|\le |G|^{-1}\left( |G(\tilde{f}(a))| +|x-a|^{1-\frac{1}{2q}}\sqrt{2q-1}\left( \int_a^b |g(\tau)|^qd\tau)\right)^{\frac{1}{q}}     \right)
\end{eqnarray*}
for every $x\in (a,b)$.\\
Analogously, if $b\in \r$, $G$ is of constant sign, $|G|$ is
 strictly increasing and $\lim_{r\to\infty}|G(r)|=\infty$ then
$f$ extends to a continuous function up to $b$ (denoted by $\tilde{f}$) and we have
\begin{eqnarray*}
|f(x)|\le |G|^{-1}\left( |G(\tilde{f}(b))| + |b-x|^{1-\frac{1}{2q}}\sqrt{2q-1}\left( \int_a^b |g(\tau)|^qd\tau)\right)^{\frac{1}{q}}     \right)
\end{eqnarray*}
whenever $x\in (a,b)$.
\item[v)] if $a\in \r$,  $G$ is monotonic and $\tau\circ G^{-1}$ is proper (i.e. it maps bounded sets into bounded sets)  then $f\in W^{2,q}(a,x)$
for any $x\in (a,b)$ and
\begin{eqnarray*}
&~&\left( \int_{a}^x |f^{''}(x)|^qdx\right)^{\frac{1}{q}} \le D_c (x) \left( \int_{a}^x |g(x)|^qdx\right)^{\frac{1}{q}}, \hbox{\rm where}\\
D_c (x)&:=& {\rm sup} \{  (\tau\circ G^{-1})(y): y\in (-A_c(x),A_c(x))\} \nonumber\\
 A_c(x)&:=& \left( |(G\circ f)(c)|+ |x-c|^{1-\frac{1}{2q}} \sqrt{2q-1}\left(\int_a^b |g(x)|^qdx\right)^{\frac{1}{2q}}\right).\nonumber
\end{eqnarray*}
Number $c\in (a,x)$ can be chosen arbitrary.

In particular if also $b\in\r$ then $f\in W^{2,q}(a,b)$  and
\begin{eqnarray*}
&~&\left( \int_{a}^b |f^{''}(x)|^qdx\right)^{\frac{1}{q}} \le D_c \left( \int_{a}^b |g(x)|^qdx\right)^{\frac{1}{q}}, \hbox{\rm where}\\
D_c&:=& {\rm sup} \{  (\tau\circ G^{-1})(y): y\in (-A_c,A_c)\} \nonumber\\
 A_c&:=& \left( |(G\circ f)(c)|+ ({\max} \{ |c-a|, |b-c|\})^{1-\frac{1}{2q}} \sqrt{2q-1}\left(\int_a^b |g(x)|^qdx\right)^{\frac{1}{2q}}\right),\nonumber
\end{eqnarray*}
number $c\in (a,b)$ can be chosen arbitrary.

\end{description}
\end{prop}

\noindent
{\bf Proof.}
{\bf i): }
By Lemma \ref{lemat} we know that  function $h$ satisfies the identity (\ref{srod}).
Therefore i) follows from Proposition \ref{winosekdodatk}.\\
{\bf ii): } By construction $G(\cdot)$ is locally Lipschitz, so that  $(G(f))^{'}=h^\frac{1}{2q}(f)f^{'}$ and by  already proven part i),  it belongs to
$L^{2q}(a,b)$.
Therefore $F=G(f)\in W^{1,2q}(r,R)$ for all $a<r<R<b$ (if $a$ and $b$ are finite numbers then $F\in W^{1,2q}(a,b)$) and it suffices to apply
 Morrey--Sobolev inequality  (\cite{ma}, Theorem 1.4.5, part (f)):
\begin{eqnarray*}
\frac{|F(x)-F(y)| }{|x-y|^{1-\frac{1}{2q}}}&\le &\frac{1}{|x-y|^{1-\frac{1}{2q}}}\int_x^y |F^{'}(x)| dx\le \frac{|x-y|^{1-\frac{1}{2q}}}{|x-y|^{1-\frac{1}{2q}}}
\left( \int_x^y |F^{'}(x)|^{2q}dx\right)^{\frac{1}{2q}}\\
& \le& \left( \int_x^y |f^{'}(x)|^{2q}h(f(x))dx \right)^{\frac{1}{2q}} \le \left( (\sqrt{2q-1})^{2q}\int_a^b |g(x)|^qdx\right)^{\frac{1}{2q}},
\end{eqnarray*}
where $x,y\in (a,b), x<y$.\\
{\bf iii):} This  is an obvious consequence of part ii).\\
{\bf iv):} We prove first statement only, as the second one is an obvious modification of the presented arguments.
 We may assume that $G$ is nonnegative (the same estimations as for $G$ hold for $|G|$).
Let us consider an arbitrary sequence $r_n\to a, r_n>a$. By the already proven part ii) we deduce that $F(r_n)$ is a Cauchy sequence, thus it is convergent.
It is easy to see that the limit does not depend on the choice of the sequence $r_n$. Therefore $F$ extends to a continuous function at $a$
(denote it by $\tilde{F}$) and
\(
\tilde{F}(a)=\lim_{r\searrow a} F(r).
\)
By the very definition $G$ maps $\r_+$ to some infinite interval $B=(d,\infty)$ and $G$ is increasing.
Therefore if $\tilde{F}(a)\in B$ then $G^{-1}$ is well defined, bounded in the neighborhood of $\tilde{F}(a)$ and
$ \lim_{r\searrow a} f(r) =\lim_{x \searrow a} G^{-1}(F(x))=G^{-1}(\tilde{F}(a))\in (0,\infty)$. In the other case $\tilde{F}(a)$ is one of the endpoints
of $B$ and it is finite, so that $\tilde{F}(a)=d\in \r$. By the monotonicity argument $G$  extends to the continuous function at $0$ and
the extension $\tilde{G}$ satisfies $\tilde{G}(0)=d$. Consequently $\lim_{x\searrow a}f(x)=\lim_{x\searrow a} G^{-1}(F(x))=\lim_{y\searrow d}G^{-1}(y)= 0$.
In both situations the assertion follows. The estimation follows now from part iii).\\
{\bf v):} We prove first assertion only.
 We have: $\tau (f) =\tau\circ G^{-1}\circ (G (f))$ and $|G(f)|$ is  no bigger than $A_c(x)$ on every interval $(a,x)$ according to part ii).
 It follows that $\tau (f)$ is bounded $(a,x)$ and  $|\tau (f)(x)|\le D_c(x)$, with an arbitrary chosen $c\in (a,x)$.
Therefore $f^{''}=g\cdot \tau (f)\in L^q(a,x)$ and
\(
\| f^{''}\|_{L^q(a,x)}\le  D_c(x)\| g\|_{L^q(a,x)} .
\)
\\
This finishes the proof of the Proposition.

\hfill $\Box$

\begin{rem}\rm
The sign of expression: $\liminf_{R\nearrow b}|f^{'}(R)|^{2q-2}f^{'}(R){H}(f(R))$ or\\
 $\limsup_{r\searrow a}|f^{'}(r)|^{2q-2}f^{'}(r){H}(f(r))$ often can be deduced as a consequence of the monotonicity of $f$ near the endpoint, see e.g. \cite{akita}
and their references for analysis of the monotonicity.
\end{rem}

\subsubsection{The special case. Problems with homogeneous nonlinearity}\label{lustration}

We are now to consider the following ODE:
\begin{equation}\label{jednorodne}
f^{''}(x)=g(x)f(x)^\alpha, \ \hbox{\rm a.e. for}\  x\in (0,b),
\end{equation}
where  $\alpha\in \r\setminus \{ 0\}$, $0<b\le\infty$,  $g\in L^q(0,b)$ and
 $f\in W^{2,1}_{loc}(0,b)$ is  positive with yet to be specified boundary conditions.\\
  Equations with homogeneous nonlinearities naturally appear in many mathematical models. Let us mention a few of them, referring mostly to  books \cite{ar,arw} for details and many more references:
\begin{description}
\item[1)] Thomas and Fermi model found independently in 1927 to determine the electrical potential in an isolated neutral atom (see  \cite{ar}, page 121, \cite{arw}, page x, and celebrated historical articles \cite{tomas}, \cite{fermi}):
\begin{eqnarray*}
\left\{
\begin{array}{c}
y^{''}(t)= t^{\frac{1}{2}}y(t)^{\frac{3}{2}} , \ t\in (0,\infty)\\
y(0)=0, \lim_{t\to\infty} y(t)=0.
\end{array}
\right.
& {\rm or} & \left\{
\begin{array}{c}
y^{''}(t)= t^{\frac{1}{2}}y(t)^{\frac{3}{2}} , \ t\in (0,1)\\
y(0)=0,  y(1)=0.
\end{array}
\right.
\end{eqnarray*}
\item[2)] The  Emden-Fowler problem (see e.g. \cite{arw}, page x) which appears in various branches of fluid dynamics:
\begin{eqnarray*}
\left\{
\begin{array}{c}
y^{''}+\lambda q(x)y^{-\gamma}=0 , \ x\in (0,1),\gamma >0\\
y(0)=y(1)=0.
\end{array}
\right.
\end{eqnarray*}
\item[3)] Problem arising in the study of deformation shape of a membrane cap (\cite{ar}, page 123):
\begin{eqnarray*}
\left\{
\begin{array}{c}
u^{''}(t)= \frac{1}{t^{3}}\left( \frac{\lambda^2}{8t^{\gamma-2}} -\frac{1}{32 u(t)^2}  +\frac{\mu}{4u(t)} \right),\ 1<t<\infty\\
a_0u(1)-a_1u^{'}(1)=A,  \ u(t)\ \hbox{\rm is bounded as}\ t\to\infty .
\end{array}
\right.
\end{eqnarray*}
If one simplifies this model by assuming $\mu=\lambda =0$, we get the ODE:\[ u^{''}(t)= -\frac{1}{32\cdot t^{3}} \frac{1}{ u(t)^2}.\]
\item[4)] Problem arising in the study of thin film of viscous fluid over a solid surface (\cite{ar}, page 297):
\begin{eqnarray*}
\left\{
\begin{array}{c}
y^{''}(u)=- \frac{2}{u^2}\frac{1}{\sqrt{y}} , \ u\in (1,\infty)\\
y(1)=0,\  \lim_{u\to\infty} y(u)=0.
\end{array}
\right.
\end{eqnarray*}
\item[5)] One parameter family of logistic equations:
\[
u^{''} +af(x)u- b(x)u^p=0,\ x\in (0,\infty)
\]
where $p>1, a\in\r, b\ge 0$ is of special interest in mathematical biology (see e.g. \cite{lixi}, \cite{afbr} and their references).
In case $a=0$ or $b\equiv 0$ it reduces to ODE (\ref{jednorodne}).
\end{description}

\noindent
We start with the following lemma.

\begin{lem} When $\tau (\lambda)=\lambda^\alpha$, $\tau: (0,\infty)\rightarrow (0,\infty)$, $q\ge 1$,
$\alpha \neq -1+\frac{1}{q}$, the transformations of $\tau$ in Definition \ref{def2} are given by: $k_{q,\tau}(\lambda)= \lambda^{\frac{\alpha q}{q-1}}$ and
  $K_{q,\tau}(\lambda)= \frac{1}{1+\frac{\alpha q}{q-1}}\lambda^{\frac{\alpha q}{q-1} +1}$ when $q>1$ (in particular $K_{q,\tau}$ is of constant sign)  and
\begin{eqnarray*}
h_{q,\tau}(\lambda) &:=&
 |q-1+\alpha q|^q \lambda^{-q(\alpha +1)}, \\
H_{q,\tau}(\lambda) &:=&
-{\rm sgn}(q-1+\alpha q) |q-1+\alpha q|^{q-1} \lambda^{-q(\alpha +1)+1},\\
G_{q,\tau}(\lambda) &:=& \left\{
\begin{array}{ccc}
 |q-1+\alpha q|^{\frac{1}{2}}\frac{2}{1-\alpha} \lambda^{\frac{1-\alpha}{2}} & {\rm if}&  \alpha\neq 1\\
 (2q-1)^{\frac{1}{2}}\ln \lambda &{\rm if}& \alpha=1
\end{array}
\right.
\end{eqnarray*}
\end{lem}

\begin{rem}\rm
In case $\alpha = -1+\frac{1}{q}$ we have $k_{q,\tau}(\lambda)= \lambda^{-1}$, so that any its locally absolutely continuous primitive
$K_{q,\tau}$ is a translation of logarithmic function. Therefore it cannot be of constant sign.
\end{rem}

\begin{rem}\rm
If $\alpha =1$ function $|G_{q,\tau}|$ is not monotonic.
\end{rem}

\noindent
Applying Proposition \ref{regul0} we arrive at the following result. Easy proof is left to the reader.

\begin{prop}\label{regu11}
Suppose that $1\le q<\infty$, $\alpha \neq -1+\frac{1}{q}$, $\kappa =-{\rm sign} (\alpha + 1-\frac{1}{q})$,  $0<b\le\infty$,  $g\in L^q(0,b)$ and
let $f\in W^{2,1}_{loc}(0,b)$ is a positive solution of the following ODE:
{\small
\begin{equation*}
\left\{
\begin{array}{c}
f^{''}(x)=g(x)(f(x))^\alpha \ \hbox{\rm a.e. on}\ (0,b)\\
 \liminf_{R\nearrow b}\kappa |f^{'}(R)|^{2q-2}f^{'}(R)(f(R))^{-q(\alpha+1)+1}
 - \limsup_{r\searrow 0} \kappa |f^{'}(r)|^{2q-2}f^{'}(r)(f(r))^{-q(\alpha+1)+1} \le 0.
 \end{array}
 \right.
\end{equation*}
}
Then we have
\begin{description}
\item[i)] $|f^{'}|^{2q}|f|^{-q(\alpha +1)} \in L^1(0,b)$ and
\begin{eqnarray*}
\int_0^b |f^{'}(x)|^{2q}|f(x)|^{-q(\alpha +1)} dx &\le& \left( \frac{(2q-1)}{|q-1+\alpha q|}\right)^q \int_0^b |g(x)|^q dx,
\end{eqnarray*}
\item[ii)]
\begin{eqnarray*}
&~&{\rm sup}\left\{ \frac{|(f(x))^{\frac{1-\alpha}{2}} - (f(y))^{\frac{1-\alpha}{2}} |}{|x-y|^{1-\frac{1}{2q}}} : x,y\in (0,b)\right\} \le A_q \left( \int_0^b |g(x)|^q dx\right)^{\frac{1}{2q}},\
\\
&~& A_q =
 \sqrt{ 2q-1} |q-1+\alpha q|^{-\frac{1}{2}}|\frac{1-\alpha}{2}|,
\end{eqnarray*}
\item[iii)] If $\alpha <1$ then $\lim_{r\searrow 0}f(r)=:f(0)$ exists and the following estimation holds:
\begin{eqnarray*}
|f(x)|&\le& \left\{ (f(0))^{\frac{1-\alpha}{2}} + A_q |x|^{1-\frac{1}{2q}} \left( \int_0^b |g(x)|^q dx\right)^{\frac{1}{q}}\right\}^{\frac{2}{1-\alpha}}.
\end{eqnarray*}
\end{description}
\end{prop}

\begin{rem}\rm Under assumptions of Proposition \ref{regu11} if $\alpha > 1$ then we necessarily have $\liminf_{r\searrow 0} f(r)>0$. This follows directly from assertion ii).
\end{rem}

\begin{rem}\rm
If $\alpha >\frac{1}{q}-1$ ($\kappa =-1$) the boundary condition reads as:
\[
- \limsup_{R\nearrow b} |f^{'}(R)|^{2q-2}f^{'}(R)(f(R))^{-q(\alpha+1)+1}
 + \liminf_{r\searrow 0} |f^{'}(r)|^{2q-2}f^{'}(r)(f(r))^{-q(\alpha+1)+1} \le 0.
\]
It is satisfied when for example $f$ is increasing near the endpoint $b$, $f(0)>0$ and
  we have either: $f$ is decreasing near $0$ or
 $f$ is $C^1$ up to $0$ and  $f^{'}(0)=0$.
\end{rem}

\noindent
Obviously, various extensions and modifications of Proposition \ref{regul0} can be applied according to the required situation. One can deal with nonnegative but not necessarily strictly positive solutions (``dead core solutions'', \cite{diaz}),
 or for example with sign changing solutions, appearing typically in  nonlinear eigenvalue problems,  \cite{lind}, \cite{annane}.
In such  cases we expect to apply Propositions: \ref{goal2}, \ref{goracs}, \ref{gorac}, \ref{doublezeroes}, \ref{goracs2}, \ref{goal111}.
We skip the precise formulations, hoping that the proposed approach will serve in the future in such cases as well.


\begin{thebibliography}{99}

\bibitem{akita} Adamowicz, T.,  Ka\l{}amajska, A., {\em Maximum principles and nonexistence results for radial solutions to
equations involving $p$-Laplacian},
Math. Methods Appl. Sci. 33, no. 13 (2010),  1618--1627.

\bibitem{afbr} Afrouzi, G.A.,  Brown, K.J., {\em  On a difusive logistic equation,} J. Math. Anal. Appl., {\bf 225} (1998), 326-339.

\bibitem{ar} Agarwal, R. O'Regan, D., {\em  Singular differential and
integral equations with applications.} Kluwer Academic
Publishers, Dordrecht, 2003.

\bibitem{arw} Agarwal, R. O'Regan, D., Wong P.,  {\em Positive solutions of Differential, Difference and Integral Equations}, Kluwer Academic
Publishers, Dordrecht, 1999.

\bibitem{annane} Anane, A., Tsouli, N.,
{\em On the second eigenvalue of the p-Laplacian,}
Nonlinear partial differential equations (F\'es, 1994), 1--9, Pitman Res. Notes Math. Ser., 343,
Longman, Harlow, 1996.


\bibitem{bloom} Bloom, S., {\em First and second order Opial inequalities,}  Studia Math.  {\bf 126}(1)  (1997),   27--50.

\bibitem{bgm} Buttazzo, G. Giaquinta, M., Hildebrandt, S.,
{\em One-dimensional variational problems. An introduction.} Oxford
Lecture Series in Mathematics and its Applications, 15. The
Clarendon Press, Oxford University Press, New York, 1998.


\bibitem{bbch} Brown, R.,  Burenkov, V., Clark, S., Hinton, D.,
{\em Second order Opial inequalities in ${\cal L}^p$ spaces and applications,}
Rassias, Themistocles M. (ed.) et al., Analytic and geometric inequalities and applications. Dordrecht: Kluwer Academic Publishers. Math. Appl., Dordr. 478, 37-52 (1999).


\bibitem{diaz} Díaz, J. I., {\em Nonlinear partial differential equations and free boundaries. Vol. I. Elliptic equations,} Research Notes in Mathematics, 106. Pitman (Advanced Publishing Program), Boston, MA, 1985.

\bibitem{fermi} Fermi, E., {\em Un methodo statistico par la determinazione di alcune properit\'a dell' atoma}, Rend. Accad. Naz. del Lincei Cl. Sci. Fis. Mat. e Nat.
{\bf 6} (1927), 602-607.

\bibitem{fitz} FitzGerald, C.H., {\em Opial-type inequalities that involve higher order derivatives,}
General inequalities 4, Mem. E. F. Beckenbach, 4th Int. Conf., Oberwolfach/Ger. 1983, ISNM 71, 25-36 (1984).


\bibitem{hmv} Hernández, J., Mancebo, F.J., Vega, J.M., {\em Nonlinear singular elliptic problems: recent results and open problems.} Nonlinear elliptic and parabolic problems, 227--242, Progr. Nonlinear Differential Equations Appl., 64, Birkhäuser, Basel, 2005.

\bibitem{g}  Gagliardo, E., {\em  Ulteriori  propriet\`a  di  alcune  classi  di
funzioni  in  pi\`u  variabili} (in Italian),  Ricerche   Mat.
{\bf 8} (1959), 24--51.

\bibitem{kp1} Ka\l amajska, A., Pietruska-Pa\l uba, K., {\em
Gagliardo-Nirenberg inequalities in Orlicz spaces}, Indiana Univ. Math. J. {\bf 55}(6) (2006), pp.
1767--1789.

\bibitem{lind}  Lindqvist P., {\em  On the equation ${\rm div}(\vert \nabla u\vert \sp{p-2}\nabla u)+\lambda \vert u\vert \sp{p-2}u=0$},
 Proc. Amer. Math. Soc. {\bf 109}, no.1 (1990),  157-164. Addendum, ibid. 116 (1992), 583-584.


\bibitem{lixi} Li Ma,  Xingwang Xu, {\em Positive solutions of a logistic equation on unbounded intervals,} Proc. Amer. Math. Soc. {\bf 130}, No. 10 (2002) , 2947--2958


\bibitem{ma}   Mazy'a,  V. G.,  {\em Sobolev Spaces,} Springer- Verlag 1985.


\bibitem{mazjakufner}   Mazy'a,  V. G., Kufner A., {\em Variations on the theme of the inequality $(f')^2\leq 2f\,  {\rm sup} | f^{''}|$,}
Manuscripta Math. {\bf 56}, no. 1 (1986), 89--104.


\bibitem{n1}   Nirenberg, L, {\em On elliptic partial differential
equations},  Ann.\ Scuola Norm.\ Sup.\ di Pisa, {\bf 13} (1959),
115-162.



\bibitem{opial} Opial, Z. {\em Sur une inégalité}, (French)  Ann. Polon. Math.  {\bf 8}  (1960),  29--32.

\bibitem{oregan} O'Regan, D. {\em Theory of Singular Boundary Problems.} World Scientific, Singapore, 1994.

\bibitem{tomas} Thomas, L.H., {\em The calculation of atomic fields,} Proc. Camb. Phil. Soc. {\bf 23} (1927), 1473-1484.



\end{thebibliography}
\end{document}